\documentclass[12pt,a4paper]{article}
\usepackage[utf8]{inputenc}
\usepackage{amsmath}
\usepackage{amsfonts}
\usepackage{amssymb}
\usepackage{hyperref}

\usepackage{geometry}
\usepackage[all]{xy}

\newcommand*{\ca}{\mathcal}
\newcommand*{\ba}{\mathrm}
\newcommand*{\arr}{\rightarrow}
\newcommand*{\map}{\mapsto}
\newcommand*{\iso}{\cong}
\newcommand*{\we}{\simeq}
\newcommand*{\qis}{\overset{\sim}{\rightarrow}}
\newcommand*{\D}{{[d]}}

\begin{document}

\title{Higher order Hochschild cohomology of schemes}
\date{Lucas DARBAS\\ \ \\ \ \\ Laboratoire de Mathématiques Jean Leray\\ \ \\ Nantes Université}
\maketitle

\noindent\textbf{Abstract} : We show that Higher Hochschild complex \cite{Pir} associated to a connected pointed simplicial set commutes with localization of commutative algebras over a field of characteristic zero. Then, we define in two ways higher order Hochschild cohomology of schemes over a field of characteristic zero. Originally, we can take the hyperext functor of the sheaf associated to Higher Hochschild presheaf. We obtain a Hodge decomposition for higher order Hochschild cohomology of smooth algebraic varieties over a field of characteristic zero which generalizes Pirashvili's Hodge decomposition. We can also define the Higher Hochschild cohomology of order $d$ of a separated scheme by taking the ext functor of its structure presheaf over the Higher Hochschild presheaf of order $d-1$. These two definitions are really close to those of Swan \cite{Swa} for classical Hochschild cohomology, but our tools are model categories and derived functors. We also generalize the equivalence of Swan's definitions to any separated schemes over a field.\\

\tableofcontents

\section*{Acknowledgements}
\
\indent This paper is a substantial part of my PhD about sheaf properties of the Higher Hochschild complex in Laboratoire Jean Leray of Nantes Université. I want to thank F. Wagemann and H. Abbaspour for supervising my PhD. I also want to thank G. Powell and G. Ginot for supporting me during my research.\\

\ \\

\section{Introduction}
\
\indent The \textit{Hochschild complex} $C(A)$ of an associative algebra $A$ is the chain complex 
\begin{center}
$\xymatrix @!0 @R=9mm @C=3.5cm {\cdots \ar[r]^-{d_0-d_1+d_2-d_3} & A\otimes A\otimes A \ar[r]^-{d_0-d_1+d_2} & A\otimes A \ar[r]^-{d_0-d_1} & A }$
\end{center}
\begin{center}
$d_i (a_0 \otimes\cdots\otimes a_n )={\left\{
    \begin{array}{ll}
        a_0 \otimes\cdots\otimes a_i a_{i+1}\otimes\cdots\otimes a_n & \mbox{if } 0\leq i<n\\
        a_na_0\otimes a_1 \otimes\cdots\otimes a_{n-1} & \mbox{if } i=n
    \end{array}
\right.}$
\end{center}
It is closely related to the simplicial circle $\mathbf{S}^1$
\begin{center}
$\xymatrix{\cdots \ar@<12pt>[r] \ar@<4pt>[r] \ar@<-4pt>[r] \ar@<-12pt>[r] & \{\ast ,s_0\sigma ,s_1\sigma\} \ar@<8pt>[r] \ar[r] \ar@<-8pt>[r] \ar@<8pt>[l] \ar[l] \ar@<-8pt>[l] & \{\ast ,\sigma\} \ar@<4pt>[r] \ar@<-4pt>[r] \ar@<4pt>[l] \ar@<-4pt>[l] & \{\ast\} \ar[l] }$
\end{center}
with in degree $n$ a copy of $A$ for each $n$-simplex. The idea of Loday \cite[6.4]{Lod} and Pirashvili \cite{Pir} was to replace the circle by any degreewise finite simplicial set $K$ when $A$ is a commutative algebra in order to build a new simplicial complex $C(K,A)$ and then to obtain a new (co)homology functor $H(K,A)$. In \hyperlink{3}{Section 3}, we will use skeleton of simplicial sets and derived tensor product \cite{GTZ2} to show that it commutes with localization in homology if the simplicial set $K$ is connected.\\

\noindent\hyperlink{3.2.2}{\textbf{Theorem 3.2.2.}} Let $K$ be a connected pointed simplicial set and $S$ a multiplicative system of an algebra $A$. The canonical morphism of CDGA
\begin{center}
$S^{-1}C(K,A)\arr C(K,S^{-1}A)$
\end{center}
is a quasi-isomorphism.\\

Pirashvili was particularly interested in the case of \textit{simplicial d-spheres} $\mathbf S^d$ and generalized the \textit{$\lambda$-decomposition} of Loday \cite[6.4.5]{Lod} to a \textit{Hodge decomposition} depending only on the parity of $d$  
\begin{center}
$H_n\big(\textbf{S}^d,A\big)\iso\bigoplus\limits_{i+jd=n}H^{(j)}_{i+j}(A)$
\end{center}
when $A$ is over a field of characteristic zero \cite[2.5]{Pir}.\\

\indent On the other hand, Swan defined in three different ways the Hochschild cohomology of a scheme $X$ \cite{Swa}. Following an idea of Grothendieck and Loday, he used the sheaf $\ca C_X$ associated to the presheaf $U\map C\big(\ca O_X(U)\big)$ and the hyperext functor $\ba{\mathbb{E}xt}_{\ca O_X}(\ca C_X,\ \ \ )$ \cite[2]{Swa}. This construction can easily be done with the simplicial complexes $C(K,A)$. In \hyperlink{4}{Section 4}, we will be interested in the cohomological functor $H_\D(X,\ \ \ )$ corresponding to the \textit{$d$-spheres} $\mathbf S^d$. We will show that there exists a \textit{Hodge decomposition for smooth varieties} and that the functor $H_\D(X,\ \ \ )$ and this decomposition are the same as Pirashvili's for affine schemes.\\

\newpage

\noindent\hyperlink{4.1.5}{\textbf{Theorem 4.1.5.}} Let $d$ be a positive integer, $X$ a separated smooth scheme over a field of characteristic zero and $\ca F$ an $\ca O_X$-module. There is a natural isomorphism
\begin{center}
$H_\D^n(X,\ca F)\iso\bigoplus\limits_{p+jd=n}H^p(X,\ca D_X^j\otimes\ca F)$
\end{center}
where $\ca D_X^j$ is the dual sheaf of $\Omega^j_X$ if $d$ is odd and of $\ba{Sym}^j_{\ca O_X} \Omega^1_X$ if $d$ is even.\\

If the ground ring is a field then the Hochschild cohomology of $A$ is given by the ext functor $\ba{Ext}_{A^e}(A,\ \ )$ \cite[1.5.]{Lod}. Guided by this fact, Gerstenhaber and Schack defined the Hochschild cohomology of a diagram of algebras $\ca O$ by the ext functor $\ba{Ext}_{\ca O^e}(\ca O,\ \ \ )$ \cite[21.0]{GS}. Swan took the structure sheaf of a separated scheme $X$ restricted to its affine open sets as a diagram of algebras to get another definition of the Hochschild cohomology of $X$ \cite[3]{Swa}. In \hyperlink{5}{Section 5}, we will adapt this definition for the higher order and we will show that it gives us the same functor $H_\D(X,\ \ \ )$ on the quasi-coherent sheaves.\\

\noindent\hyperlink{5.2.3}{\textbf{Theorem 5.2.3.}} Let $d$ be a positive integer, $X$ a separated scheme over a field of characteristic zero with structure presheaf $\ca O$ and $\ca F$ a quasi-coherent $\ca O_X$-module thought of as a preasheaf on $X$. There is a natural isomorphism
\begin{center}
$H^n_\D(X,\ca F)\iso H^n\big(\mathbb{R}\ba{Hom}_{C^{[d-1]}(\ca O)}(\ca O,\ca F)\big)$
\end{center}
where $C^{[d-1]}(\ca O)$ is the presheaf of CDGA on $X$ given by $U\map C\big(\mathbf{S}^{d-1},\ca O_X(U)\big)$.\\

\indent This same remark led Swan to give a more simple definition of the \textit{Hochschild cohomology} of a scheme $X$ with diagonal morphism $\delta :X\arr X\times X$ by taking the ext functor $\ba{Ext}_{\ca O_{X\times X}}(\delta_\ast\ca O_X,\delta_\ast\ \ \ )$ \cite[1]{Swa}. He showed that this functor is the same as Grothendieck-Loday's for quasi-projective schemes over a field. In \hyperlink{6}{Section 6}, we will generalize this fact to any separated scheme over a field.\\

\noindent\hyperlink{6.1.4}{\textbf{Theorem 6.1.4.}} Let $X$ be a separated scheme over a field with diagonal morphism $\delta:X\arr X\times X$ and $\ca F$ an $\ca O_X$-module. There is a natural isomorphism
\begin{center}
$\ba{\mathbb{E}xt}^n_{\ca O_X}(\ca C_X,\ca F)\iso \ba{Ext}^n_{\ca O_{X\times X}}\big(\delta_*\ca O_X,\delta_*\ca F\big) $
\end{center}
\ \\
\section{Notations and conventions}
\
\noindent 1. \textit{Homological algebra.} Our complexes are homologicaly graded, but we can talk about the cohomology of a complex $C$ with the convention $H^n(C)=H_{-n}(C)$. We denote by $\ba{Ch}(\ca A)$ the category of complexes in an abelian category $\ca A$ and $\ba{Hom}_{\ba{Ch}(\ca A)}$ the \textit{hom complex} in $\ba{Ch}(\ca A)$ \cite[5.1]{AX}. If $\ca A$ is a Grothendieck category \cite[1]{Hov2} then $\ba{Ch}(\ca A)$ has an \textit{injective} model structure which is \textit{cofibrantly generated}, \textit{proper} and such that the weak equivalences are the quasi-isomorphisms and the cofibrations are the monomorphisms \cite[2.2]{Hov2}. We denote by $\ba{Vect}$ the category of \textit{vector spaces} over a base field $k$ of caracteristic zero, $\otimes={\otimes}_k$, $\ba{Ch}$=$\ba{Ch(\ba{Vect})}$ and we call complexes in $\ba{Vect}$ \textit{chain complexes}.\\

\noindent 2. \textit{Simplicial objects.} Let $\mathbf{\Delta}$ be the category with objects the finite sets $[n]=\{0,\dots,n\}$ for all $n\in\mathbb{N}$ and with morphisms the non-decreasing applications. A simplicial object of a category $\ca C$ is a functor $\mathbf{\Delta}^\ba{op}\arr\ca C$ and a cosimplicial object of $\ca C$ is a functor $\mathbf{\Delta}\arr\ca C$. The category of simplicial objects of $\ca C$ is denoted by $\ba s\ca C$. A \textit{finite} simplicial set is a simplicial object of the category of finite sets $\ba{Fin}$. We denote by $\ba{Fin'}$ the category of fnite pointed sets. For any $d\in\mathbb N$, let $\mathbf{\Delta}^d=\ba{Hom}_\mathbf{\Delta}\big(\ \ \ ,[d]\big)$ the \textit{standard $d$-simplex}, $\partial\mathbf{\Delta}^d\subset\mathbf{\Delta}^d$ the simplicial subset generated by the non-degenerated simpleces $\sigma\neq 1_\D$ its \textit{border} and $\ast=\mathbf{\Delta}^0$ the \textit{point}.
We define the \textit{$d$-sphere} by the simplicial pushout $\mathbf{S}^d=\ast\cup_{\partial\mathbf{\Delta}^d}\mathbf{\Delta}^d$. The category $\ba{sSet}$ of simplicial sets has a \textit{Kan} model structure which is cofibrantly generated, proper and such that the cofibrations are injections \cite[3.6.5, 3.2.2]{Hov1} \cite[13.1.13]{Hir}.\\

\noindent 3. \textit{CDGA.} A monoid in a monoidal category $(\ca M,\otimes ,1)$ is an object $A$  with two morphisms $A\otimes A\arr A$ and $1\arr A$ and an $A$-module is an objet $M$ of $\ca M$ with a morphism $A\otimes M\arr M$. They have to satisfy usual diagram relations \cite[VII, 3, 4]{ML2}. We denote by $A$-mod the category of $A$-modules, $\ba{Sym}^n_A$ the $n^\ba{th}$ symmetric power over $A$ and $\bigwedge^n_A$ the $n^\ba{th}$ exterior power over $A$ for any $n\in\mathbb N$. We call commutative monoids in $(\ba{Vect},\otimes ,k)$ \textit{algebras} and we denote by $\ba{Alg}$ their category. For any algebra $A$, we denote by $\ba{Spec}(A)$ the set of prime ideals of $A$, $\Omega^1_A$ the $A$-module of differentials of $A$ \cite[III, 10.11]{A} and $\Omega^n_A=\bigwedge^n_A\Omega^1_A$ for any $n\in\mathbb N$. For any multiplicative system $S$ of $A$, we denote by $S^{-1}A$ the fraction ring of $A$ defined by $S$ \cite[II, 1.1]{AC} and $S^{-1}M=S^{-1}A\otimes_A M$ for any $A$-module $M$. If $S=\{1,s,s^2,\dots ,s^n,\dots \}$ we use the notations $A_s=S^{-1}A$ and $M_s=S^{-1}M$. We call commutative monoids in $(\ba{Ch},\otimes , k)$ \textit{CDGA} and we denote by $\ba{CDGA}$ their category. Since $k$ is a field of caracteristic zero, it is a cofibrantly generated model category such that the weak equivalences are the quasi-isomorphisms and the fibrations are the surjections \cite[4.1.1, 4.2.5.2]{Hin}. Any simplicial algebra $A:\mathbf{\Delta}^\ba{op}\arr\ba{Alg}$ can be thought as a CDGA $A=(A_n)_{n\in\mathbb N}$ with differential $d_n=\sum^n_{i=0}(-1)^id_i$ \cite[8.2.1]{Wei} and multiplication the composition of the \textit{suffle product} \cite[VIII, 8.8]{ML1} and the multiplications of $A$.\\

\noindent 4. \textit{Sheaves.} All the sheaves on a topological space $X$ take values in $\ba{Vect}$. For any sheaf $\ca F$ on $X$, we use notation $\Gamma\ca F=\ca F(X)$. A ringed space $X$ is a topological space with a sheaf $\ca O_X$ with values in $\ba{Alg}$. We denote by $\ca Hom_{\ca O_X}$ the \textit{hom sheaf} \cite[4.1]{Toh}, $\ca Ext_{\ca O_X}$ the \textit{ext sheaf} \cite[4.2]{Toh} and $\ca F^\vee=\ca Hom_{\ca O_X}(\ca F,\ca O_X)$ the \textit{dual sheaf} of an $\ca O_X$-module $\ca F$. If it is not confusing then we can use the notation $\otimes ={\otimes}_{\ca O_X}$. By convention, all schemes are over $\ba{Spec}(k)$. We denote by $\times$ the product of schemes over $\ba{Spec}(k)$ \cite[3.2.1]{EGA}, $\Omega^1_X$ the $\ca O_X$-module of $1$-differentials of $X$ over $\ba{Spec}(k)$ \cite[16.3.1]{EGAIV4} and $\Omega^n_X=\bigwedge^n_{\ca O_X}\Omega^1_X$ for any $n\in\mathbb N$.\\

\hypertarget{3}{\section{Higher Hochschild functor}}
\
\indent In all of this section, if nothing else is specified, $k$ is a field of caracteristic zero. Recall that algebras are associatives and commutatives over $k$ with identity.\\

Let $A$ be an algebra. The Loday functor \cite[6.4.4]{Lod}
\begin{center}
$\ca L(A,A):\ba{Fin}\arr \ba{Alg}$
\end{center}
associates to each finite set $K$ the tensor product $A^{\otimes K}$ and to each map $f:K\arr L$ the morphism of algebras $f_*:A^{\otimes K}\arr A^{\otimes L}$ given by
\begin{center}
$f_*\Big(\underset{x\in K}{\bigotimes}a_x\Big)=\underset{y\in L}{\bigotimes}\Big(\underset{f(x)=y}{\prod}a_x\Big)$
\end{center}
The idea of Pirashvili \cite[2.1]{Pir} is to build a simplicial algebra $C(K,A)$ associated to any finite simplicial set $K$ by taking the composition
\begin{center}
$C(K,A):\xymatrix @!0 @R=9mm @C=2.5cm { \mathbf{\Delta}^\ba{op} \ar[r]^{K} & \ba{Fin} \ar[r]^{\ca L(A,A)} & \ba{Alg} }$
\end{center}
which induices a CDGA with homology denoted by $H(K,A)$. Remark that if $A$ is a CDGA then one can define $C(K,A)$ by taking the total complex of $\big(A^{\otimes K_p}\big)_q$ \cite[3.1]{GTZ2}. We can also define $C(K,A)$ if $K$ is not finite by taking the colimit of $C(L,A)$ for all finite simplicial subsets $L\subset K$. Thus, we defined the \textit{Higher Hochschild functor}
\begin{center}
$C : \ba{sSet}\times\ba{CDGA}\arr\ba{CDGA}$
\end{center}
It preserves weak equivalences and homotopy colimits of simplicial sets \cite[3.2]{GTZ2}.\\

\noindent\hypertarget{3.0.1}{\textbf{Example 3.0.1.}} Let $A$ be a CDGA. For any diagram of simplicial set
\begin{center}
$X\leftarrow Z\arr Y$
\end{center}
ve have a natural weak equivalence of CDGA 
\begin{center}
$C(X\cup_Z^h Y,A)\we C(X,A)\overset{\mathbb{L}}{\underset{C(Z,A)}{\otimes}}C(Y,A)$
\end{center}
\ \\
\indent Because we want to work with schemes, we will only studie $C(K,A)$ with CDGA concentrated in degree 0, i.e. algebras. The following example is important to give a $A$-module structure on $C(K,A)$.\\

\noindent\hypertarget{3.0.2}{\textbf{Exemple 3.0.2.}} Let $A$ be an algebra. The chain complex $C(\ast ,A)$ is given by
\begin{center}
$\xymatrix{\cdots \ar[r]^0 & A \ar[r]^1 & A \ar[r]^0 & A }$
\end{center}
so $C(\ast ,A)$ deformation retract of $A$ as a CDGA. In particulary, if $K$ is a pointed simplicial set then $C(K,A)$ is a CDGA over $A$ :
\begin{center}
$A\arr C(\ast,A)\arr C(K,A)$ 
\end{center}
\ \\
\subsection{Higher order Hochschild Homology}
\
\indent We first introduce the \textit{higher order Hochschild complex} and its homology.\\

\noindent\hypertarget{3.1.1}{\textbf{Definition 3.1.1.}} Let $d$ be a positive integer and $A$ an algebra. We define
\begin{center}
$C^\D (A)=C(\mathbf{S}^d,A)$
\end{center}
the \textit{Higher Hochschild complex of order $d$ of $A$} and
\begin{center}
$H^\D (A)=H(\mathbf{S}^d,A)$
\end{center}
the \textit{Higher Hochschild homology of order $d$ of $A$}.\\

We will use the notations $C(A)=C^{[1]}(A)$ and $HH(A)=H^{[1]}(A)$ for the Hochschild complex of $A$ and its homology.\\

One can compute the first higher order Hochschild homology groups.\\

\noindent\hypertarget{3.1.2}{\textbf{Example 3.1.2.}} For any positive integer $d$, there is a natural isomorphism
\begin{center}
$H^\D_q(A)\iso{\left\{
    \begin{array}{ll}
        A & \mbox{if } q=0\\
        0 & \mbox{if } 0<q<d\\
        \Omega^1_A & \mbox{if } q=d
    \end{array}
\right.}$
\end{center}
This can be proved using the \textit{Hodge decomposition} of Pirashvili \cite[2.5, 1.15]{Pir} in characateristic zero, but also by hand in \textbf{any characteristic}. The computation depends on the parity of $d$ but the result is the same. Recall our simplicial model $\mathbf{S}^d$ for the sphere
\begin{center}
$\dots$\ ,\ \ $\textbf{S}^d_{d+1}=\{\ast ,s_0\sigma ,\dots , s_d\sigma\}$\ ,\ \ $\textbf{S}^d_d=\{\ast ,\sigma\}$\ ,\ \ $\textbf{S}^d_{d-1}=\{\ast\}$\ ,\ \ $\dots$\ ,\ \ $\textbf{S}^d_0=\{\ast\}$
\end{center}
To describe the differential, just look at the fibers over $\{\sigma\}$ by the face maps $d_i$ :
\begin{center}
${d_i}^{-1}\{\sigma\}={\left\{
    \begin{array}{ll}
        \{s_0\sigma\} & \mbox{if } i=0\\
        \{s_{i-1}\sigma ,s_i\sigma\} & \mbox{if } 0<i\leq d\\
        \{s_d\sigma\} & \mbox{if } i=d+1
    \end{array}
\right.}$
\end{center}
\underline{$d$ odd} : The chain complex $C^\D (A)$ is given by
\begin{center}
$\xymatrix{\cdots \ar[r] & A\otimes A^{\otimes d+1} \ar[r]^-d & A\otimes A \ar[r]^-0 & A \ar[r]^-1 & \cdots \ar[r]^-1 & A \ar[r]^-0 & A }$
\end{center}
so we just have to show that $H^\D_d(A)\iso\Omega^1_A$. Here we can use the \textit{Leibniz rule} definition of $\Omega^1_A$ \cite[1.1.9]{Lod}. The boundaries of degree $d$ are generated in $A\otimes A$ by \textit{Leibniz rule} elements of the form $1\otimes xy-y\otimes x-x\otimes y$ :
\begin{center}
$d(1\otimes x_0\otimes\cdots\otimes x_d)=\Big(\prod\limits_{j\neq 0}x_j\Big)\otimes x_0+\sum\limits_{i=1}^d (-1)^i\Big(\prod\limits_{j\neq i-1,i}x_j\Big)\otimes x_{i-1}x_i+\Big(\prod\limits_{j\neq d}x_j\Big)\otimes x_d$\\
$=\sum\limits_{i=1}^d (-1)^{i-1}\Big(\prod\limits_{j\neq i-1,i}x_j\Big)(x_i\otimes x_{i-1}-1\otimes x_{i-1}x_i+x_{i-1}\otimes x_i)$
\end{center}
Conversely, the \textit{Leibniz rule} elements of the form $1\otimes xy-x\otimes y-x\otimes y$ are boundaries :
\begin{center}
$d(1\otimes x\otimes y\otimes 1\otimes\cdots\otimes 1)=y\otimes x-1\otimes xy+x\otimes y+\sum\limits_{i=3}^{d+1}(-1)^ixy\otimes 1$\\
$=y\otimes x-1\otimes xy+x\otimes y$
\end{center}
\underline{$d$ even} : The chain complex $C^\D (A)$ is given by
\begin{center}
$\xymatrix{\cdots \ar[r] & A\otimes A^{\otimes d+1} \ar[r]^-d & A\otimes A \ar[r]^-\mu & A \ar[r]^-0 & \cdots \ar[r]^-1 & A \ar[r]^-0 & A }$
\end{center}
where $\mu$ is the multiplication of $A$, so we just have to show that $H^\D_d(A)\iso\Omega^1_A$. Here we can use the isomorphism $\Omega^1_A\iso\ba{ker}(\mu )/\ba{ker}(\mu )^2$ \cite[1.3.7]{Lod}. Since $\ba{ker}(\mu )$ is generated by the elements of the form $1\otimes x-x\otimes 1$, the boundaries of degree $d$ in $A\otimes A$ are in $\ba{ker}(\mu)^2$ :
\begin{center}
$d(1\otimes x_0\otimes\cdots\otimes x_d)=\Big(\prod\limits_{j\neq 0}x_j\Big)\otimes x_0+\sum\limits_{i=1}^d (-1)^i\Big(\prod\limits_{j\neq i-1,i}x_j\Big)\otimes x_{i-1}x_i-\Big(\prod\limits_{j\neq d}x_j\Big)\otimes x_d$\\
$=\sum\limits_{i=1}^d (-1)^{i-1}\Big(\prod\limits_{j\neq i-1,i}x_j\Big)(x_i\otimes x_{i-1}-1\otimes x_{i-1}x_i+x_{i-1}\otimes x_i)$\\
$=\sum\limits_{i=1}^d (-1)^{i-1}\Big(\prod\limits_{j\neq i-1,i}x_j\Big)(1\otimes x_{i-1}-x_{i-1}\otimes 1)(x_i\otimes 1-1\otimes x_i)$
\end{center}
\begin{flushleft}
since\hspace*{.5cm}$\sum\limits_{i=1}^d (-1)^{i-1}\Big(\prod\limits_{j\neq i-1,i}x_j\Big)(x_{i-1}x_i\otimes 1)=\Big(\prod\limits_{j=0}^{d}x_j\Big)\sum\limits_{i=1}^d (-1)^{i-1}(1\otimes 1)=0$ .
\end{flushleft}
Conversely, $\ba{ker}(\mu)^2$ is generated by boundaries :
\begin{center}
$d(1\otimes x\otimes y\otimes 1\otimes\cdots\otimes 1)=y\otimes x-1\otimes xy+x\otimes y+\sum\limits_{i=3}^{d+1}(-1)^ixy\otimes 1$\\
$=(1\otimes x-x\otimes 1)(y\otimes 1-1\otimes y)$
\end{center}
\ \\
\indent An important consequence of this computation is that any $A$-module can be thought of as a $C^\D (A)$-module using the morphism of CDGA
\begin{center}
$C^\D (A)\arr A$
\end{center}
given by the natural projection of a non-negatively graded CDGA over its homology of degree $0$.\\

To compute all degrees of the higher order Hochschild homology groups, one can use with more hypotheses the following theorem of Pirashvilli which can be found in \cite{Pir}.\\

\noindent\hypertarget{3.1.3}{\textbf{Pirashvili's Theorem 3.1.3.}} Let $d$ be a positive integer and $A$ a smooth algebra. There is a natural isomorphism
\begin{center}
$H^\D_q(A)\iso {\left\{
    \begin{array}{ll}
        \Omega^j_A & \mbox{if } q=jd\ \& \ d \mbox{ is odd}\\
        \ba{Sym}^j_A \Omega^1_A & \mbox{if } q=jd\ \& \ d \mbox{ is even}\\
        0 & \mbox{if } d\nmid q
    \end{array}
\right.}$
\end{center}
\ \\

\subsection{Localization}
\
\indent If $R$ is a commutative gradued algebra (CGA for short) and $S$ is a multiplicative system of the algebra $R_0\subset R$ then we define the localization $S^{-1}R=(S^{-1}R_0)\otimes_{R_0}R$ which is also a CGA. If $R$ is in fact a CDGA, then $S^{-1}R$ naturally has a differential. It gives rise to an exact functor
\begin{center}
$R\ba{-mod}\arr (S^{-1}R)\ba{-mod}$
\end{center}
given by the flat base change $(S^{-1}R)\otimes_R$ which we simply denote by $S^{-1}$.\\

It is well known that Hochschild homology commutes with localization. Here we give the same result about Higher Hochschild homology associated to connected pointed simplicial set.\\

\noindent\textbf{Examples 3.2.1.} \hypertarget{3.2.1.1}{1.} \textit{The circle}. By definition, we have $\mathbf S^1\we \ast\cup^h_{\partial\mathbf\Delta^1}\mathbf\Delta^1\we\ast\cup^h_{\ast\sqcup\ast}\ast$ so \hyperlink{3.0.1}{Exemples 3.0.1} and \hyperlink{3.0.2}{3.0.2} show that the Hochschild complex of an algebra $A$ is
\begin{center}
$C(A)\we A\underset{A^e}{\overset{\mathbb L}{\otimes}}A$
\end{center}
with $A^e=A\otimes A$. This well known weak equivalence exists in \textbf{any caracteristic} \cite[1.1.13]{Lod}. If $S$ is a multiplicative system of $A$ then we have $S^{-1}A\iso(S^{-1}A)^e\otimes_{A^e}A$ and so
\begin{center}
$S^{-1}C(A)\we (S^{-1}A)\underset{A^e}{\overset{\mathbb L}{\otimes}}A\we (S^{-1}A)\underset{(S^{-1}A)^e}{\overset{\mathbb L}{\otimes}}(S^{-1}A)\we C(S^{-1}A)$
\end{center}
since localization is a flat base change.\\
\hypertarget{3.2.1.2}{2.} \textit{The standard simplices}. Let $n\in\mathbb N$. Any vertex of $\mathbf\Delta^n$ gives rise to a weak equivalence of simplicial set $\ast\qis\mathbf\Delta^n$. We have $A\qis C(\mathbf\Delta^n,A)$ for any algebra $A$ by \hyperlink{3.0.2}{Exemple 3.0.2} \cite[2.4]{Pir}. Thus, if $S$ is a multiplicative system of $A$ then the commutative triangle
\begin{center}
$\xymatrix{ & S^{-1}A \ar[ld] \ar[rd] & \\
S^{-1}C(\mathbf\Delta^n,A) \ar[rr] & & C(\mathbf\Delta^n,S^{-1}A)}$
\end{center}
and \textit{Propertie 2-out-of-3} \cite[1.1.3]{Hov1} shows that $S^{-1}C(\mathbf\Delta^n,A)\we C(\mathbf\Delta^n,S^{-1}A)$.\\

\newpage

\noindent\hypertarget{3.2.1.3}{3.} \textit{Two points}. For any algebra $A$, we have $C(\ast\sqcup\ast,A)\iso C(\ast,A)^e$ \cite[3.1]{GTZ2} so by \hyperlink{3.0.2}{Exemple 3.0.2}, we have
\begin{center}
$S^{-1}C(\ast\sqcup\ast,A)\we S^{-1}A\otimes A\ncong (S^{-1}A)^e\we C(\ast\sqcup\ast,S^{-1}A)$
\end{center}
\cite[3.1]{GTZ2}.\\

\noindent\hypertarget{3.2.2}{\textbf{Theorem 3.2.2.}} Let $K$ be a connected pointed simplicial set and $S$ a multiplicative system of an algebra $A$. The canonical morphism of CDGA
\begin{center}
$S^{-1}C(K,A)\arr C(K,S^{-1}A)$
\end{center}
is a quasi-isomorphism.\\

\noindent\textit{Proof} Consider the skeleton of $K$
\begin{center}
$K^0\subset K^1\subset\cdots\subset K^n\subset K^{n+1}\subset\cdots\subset K$
\end{center}
with $K^n\subset K$ the simplicial subset generated by the non-degenerate simplices $\sigma\in K_p$ for $p\leq n$. We have two weak equivalences of simplicial sets
\begin{center}
$K^n\we K^{n-1}\cup^h_{(\partial\mathbf\Delta^n)^{\sqcup\Sigma_n}}(\mathbf\Delta^n)^{\sqcup\Sigma_n}$\hspace*{1cm}$\&$\hspace*{1cm}$K\we\underset{n\in\mathbb N}{\ba{hocolim}}\ K^n$
\end{center}
with $\Sigma_n$ the set of non-degenerated $n$-simplices of $K$ \cite[5.1.3]{Hov1}. Since localization and Higher Hochschild functor preserve homotopy colimits \cite[3.2]{GTZ2}, there is a factorization
\begin{center}
$S^{-1}C(K,A)\we\underset{n\in\mathbb N}{\ba{hocolim}}\ S^{-1}C(K^n,A)\arr\underset{n\in\mathbb N}{\ba{hocolim}}\ C(K^n,S^{-1}A)\we C(K,S^{-1}A)$
\end{center}
Thus, we just have to show that the canonical morphism of CDGA
\begin{center}
$S^{-1}C(K^n,A)\arr C(K^n,S^{-1}A)$
\end{center}
is a quasi-isomorphism from a certain rank. We work by induction on $n\geq 1$. For $n=1$, the \textit{geometric realization} \cite[3.1]{Hov1} of the simplicial set $K^1$ is a connected \textit{graph} \cite[0.1]{Hat}, so there exist a subset $E\subset\Sigma_1$ and a weak equivalence of simplicial set $K^1\we (\mathbf S^1)^{\cup^h_\ast E}$ \cite[0.7]{Hat}. By \hyperlink{3.2.1.1}{Exemple 3.2.1.1}, we have weak equivalences of CDGA
\begin{center}
$S^{-1}C(K^1,A)\we S^{-1}\Big(C(A)^{\otimes^\mathbb{L}_A E}\Big)$
\end{center}
\begin{center}
$\we\big(S^{-1}C(A)\big)^{\otimes^\mathbb{L}_{S^{-1}A} E}$
\end{center}
\begin{center}
$\qis\big(C(S^{-1}A)\big)^{\otimes^\mathbb{L}_{S^{-1}A} E}\we C(K^1,S^{-1}A)$
\end{center}
Assume that $n>1$ and choose a vertex $\ast\arr\partial\mathbf\Delta^n$. By \hyperlink{3.0.1}{Exemple 3.0.1} and induction, we have weak equivalences of CDGA
\begin{center}
$S^{-1}C(K^n,A)\we S^{-1}C(K^{n-1},A)\overset{\mathbb{L}}{\underset{C(\partial\mathbf\Delta^n,A)^{\otimes\Sigma_n}}{\otimes}} C(\mathbf\Delta^n,A)^{\otimes\Sigma_n}$
\end{center}
\begin{center}
$\qis C(K^{n-1},S^{-1}A)\overset{\mathbb{L}}{\underset{C(\partial\mathbf\Delta^n,A)^{\otimes\Sigma_n}}{\otimes}} C(\mathbf\Delta^n,A)^{\otimes\Sigma_n}$
\end{center}
\begin{center}
$\we C(K^{n-1},S^{-1}A)\overset{\mathbb{L}}{\underset{(S^{-1}C(\partial\mathbf\Delta^n,A))^{\otimes\Sigma_n}}{\otimes}} \big(S^{-1}C(\mathbf\Delta^n,A)\big)^{\otimes\Sigma_n}$
\end{center}

\newpage

\noindent\cite[3.1]{GTZ2}. Here we used the flat base change
\begin{center}
$(A^{\otimes n+1})^{\otimes\Sigma_n}\arr (S^{-1}A\otimes A^{\otimes n})^{\otimes\Sigma_n}$
\end{center}
and the commutative triangle of CDGA
\begin{center}
$\xymatrix{C(\partial\mathbf\Delta^n,A)^{\otimes\Sigma_n} \ar[r] \ar[d] & C(K^{n-1},S^{-1}A) \\
\big(S^{-1}C(\partial\mathbf\Delta^n,A)\big)^{\otimes\Sigma_n}  \ar@{-->}[ru] & }$
\end{center}
This computation show that it is enough to show that \hyperlink{3.2.2}{Theorem 3.2.2} is true for the simplicial sets $\mathbf\Delta^d$ and $\partial\mathbf\Delta^d$ for any integer $d\geq 2$ \cite[3.3.2]{Hin}. For $\mathbf\Delta^d$, it is \hyperlink{3.2.1.2}{Example 3.2.1.2}. For $\partial\mathbf\Delta^d$, we can work by induction on $d\geq 2$ using the homotopy pushout $\partial\mathbf\Delta^d\we\mathbf\Lambda^{d,0}\cup^h_{\partial\mathbf\Delta^{d-1}}\mathbf\Delta^{d-1}\we\ast\cup^h_{\partial\mathbf\Delta^{d-1}}\ast$ \cite[3.1]{Hov1}. If $d=2$ then $\partial\mathbf\Delta^2=\ast\cup^h_{\ast\sqcup\ast}\ast$ as in \hyperlink{3.2.1.1}{Example 3.2.1.1}. Assume that $d\geq 2$ and choose a vertex $\ast\arr\partial\mathbf\Delta^{d-1}$. Since localization is a flat base change, we have by \hyperlink{3.0.1}{Example 3.0.1} and induction weak equivalences of CDGA
\begin{center}
$S^{-1}C(\partial\mathbf\Delta^d,A)\we S^{-1}A\underset{C(\partial\mathbf\Delta^{d-1},A)}{\overset{\mathbb L}{\otimes}}A$
\end{center}
\begin{center}
$\we S^{-1}A\underset{S^{-1}C(\partial\mathbf\Delta^{d-1},A)}{\overset{\mathbb L}{\otimes}}S^{-1}A$
\end{center}
\begin{center}
$\we S^{-1}A\underset{C(\partial\mathbf\Delta^{d-1},S^{-1}A)}{\overset{\mathbb L}{\otimes}}S^{-1}A\we C(\partial\mathbf\Delta^d,S^{-1}A)$
\end{center}
This ends the proof.\\

\noindent\hypertarget{3.2.3}{\textbf{Corollary 3.2.3.}} Le $d$ a positive integer and $S$ a multiplicative system of an algebra $A$. The canonical morphism of CDGA
\begin{center}
$S^{-1}C^\D (A)\arr C^\D (S^{-1}A)$
\end{center}
is a quasi-isomorphism.\\

\noindent\hypertarget{3.2.4}{\textbf{Remark 3.2.4}} There is an other way to proove \hyperlink{3.2.2}{Theorem 3.2.2}. Any pointed simplicial set is weak equivalent to a \textit{minimal} pointed simplicial set \cite[9]{May}. Any connected minimal simplicial set $K$ has a single vertex \cite[9.2]{May} and one can show that $K$ is weak equivalent to the homotopy colimit of its pointed simplices $\ast\arr\mathbf\Delta^d\arr K$. Since localisation and Higher Hochschild complex preserve homotopy colimits, \hyperlink{3.2.1.2}{Example 3.2.1.2} leads to the general case.\\

\hypertarget{4}{\section{Grothendieck-Loday type definition}}
\
\indent In all of this section, if nothing else is specified, $k$ is a field of caracteristic zero. Recall that algebras are associatives and commutatives over $k$ with identity.\\

Let $X$ be a ringed space. The \textit{injective model structure} on the category $\ba{Ch}(\ca O_X)$ of complexes of $\ca O_X$-modules is a model structure such that the weak equivalences are the quasi-isomorphisms and the cofibrations are the monomorphisms \cite[2]{Hov2}. The fibrant objects of $\ba{Ch}(\ca O_X)$ are \textit{injective in the sense of Spaltenstein} \cite[1.5]{Spa} \cite[1.2.10]{Hov1} so the hom complex functor
\begin{center}
$\ba{Hom}_{\ba{Ch}(\ca O_X)}:\ba{Ch}(\ca O_X)^\ba{{op}}\times\ba{Ch}(\ca O_X)\arr\ba{Ch}$
\end{center}
has a right derived functor.\\

\newpage

Let $K$ be a pointed simplicial set. We denote by $\ca C(K,X)$ the sheaf of CDGA on $X$ associated to the presheaf of CGDA on $X$
\begin{center}
$U\map C\big (K,\ca O_X (U)\big )$
\end{center}
One can define the \textit{$K$-cohomology of $X$ with coefficients in any complex of $\ca O_X$-modules $\ca D$}, let $H(K,X,\ca D)$, by taking the cohomology of $\mathbb{R}\ba{Hom}_{\ba{Ch}(\ca O_X)}\big(\ca C(K,X) ,\ca D\big)$. In this paper, we are interested in the case where $K=\mathbf{S}^d$ for any positive integer $d$ and $\ca D$ is an $\ca O_X$-module thought of as a complex concentrated in degree $0$.\\

\noindent\hypertarget{4.0.1}{\textbf{Definition 4.0.1.}} Let $d$ a positive integer and $X$ a ringed space. We define
\begin{center}
$\ca C^\D_X=\ca C(\mathbf{S}^d,X)$
\end{center}
the \textit{Higher Hochschild sheaf of order $d$ of $X$}. Let $\ca F$ be an $\ca O_X$-module. We define
\begin{center}
$H_\D (X,\ca F)=H\Big(\mathbb{R}\ba{Hom}_{\ba{Ch}(\ca O_X)}\big(\ca C_X^\D ,\ca F\big)\Big)$
\end{center}
the \textit{Higher Hochschild cohomology of order $d$ of $X$ with coefficients in $\ca F$}.\\

We will use the notations $\ca C_X=\ca C^{[1]}_X$ and $HH(X,\ca F)=H_{[1]}(X,\ca F)$ for the Hochschild sheaf of $X$ and for the Hochschild cohomology of $X$ with coefficients in $\ca F$. When $X$ is a scheme, we recover the \textit{Grothendieck-Loday} type definition \cite[2]{Swa}. In fact, if we choose an injective resolution of $\ca O_X$-modules $\ca F\arr \ca I$, then we get a fibrant replacement $\ca I$ of the complex of $\ca O_X$-modules $\ca F$, since $\ca I$ is bounded above \cite[2.12]{Hov2}. Hence the chain complex $\ba{Hom}_{\ba{Ch}(\ca O_X)}(\ca C_X,\ca I)$ computes both the derived hom complex $\mathbb{R}\ba{Hom}_{\ba{Ch}(\ca O_X)}\big(\ca C_X,\ca F\big)$ and the hyperext functor $\mathbb{E}\ba{xt}_{\ca O_X}(\ca C_X,\ca F)$.\\

\noindent\hypertarget{4.0.2}{\textbf{Example 4.0.2.}} \textit{A non-separated scheme.} Let $A$ be an algebra over a \textbf{field of any caracteristic} and $s\in A$ an element such that the canonical morphism $A\arr A_s$ is not surjective. We define $S=\ba{Spec}(A)$, $U=D(s)=\{\mathfrak p\in S : s\notin\mathfrak p \}$ and $X=S\cup_US$ the \textit{affine scheme $S$ with doubled subvariety $Z=V(s)=S\setminus U$} \cite[I, 2.3.2]{EGA}. By assumption, $X$ is not separated \cite[I, 5.5.6]{EGA} and there is a diagram of schemes
\begin{center}
$1_S:\xymatrix{S \ar@<3pt>[r]^i  \ar@<-3pt>[r]_j & X \ar[r]^p & S}$
\end{center}
since $\Gamma\ca O_X=\Gamma\ca O_S\times_{\ca O_S(U)}\Gamma\ca O_S=A\times_{A_s}A$. The relation $i^*p^*=j^*p^*=1_{\ca O_S-\ba{mod}}$ shows that the functor $p^*$ is exact and so the pair
\begin{center}
$p^*:\xymatrix @!0 @R=5mm @C=25mm {
    \ba{Ch}(\ca O_S) \ar@<3pt>[r] & \ba{Ch}(\ca O_X) \ar@<3pt>[l]}:p_*$
\end{center}
is a Quillen adjonction \cite[2.13]{Hov2}. Moreover, for any simplicial set $K$, we have weak equivalences of sheaves of CDGA on $X$
\begin{center}
$\ca C(K,X)\iso p^*\ca C(K,S)\we\mathbb Lp^*\ca C(K,X)$
\end{center}
Deriving the enriched adjonction $\ba{Hom}_{\ba{Ch}(\ca O_S)}\circ(p^*\times 1)\iso \ba{Hom}_{\ba{Ch}(\ca O_X)}\circ(1\times p_*)$ \cite[1.3.7]{Hov1}, we obtain a natural isomorphism
\begin{center}
$H^n(K,X,\ca D)\iso H^n(K,S,\mathbb Rp_*\ca D)$
\end{center}
for any complex of $\ca O_X$-modules $\ca D$. In \hyperlink{6.1.5}{Section 6}, we will compute with more precision the cohomology groups $HH^n(X,\ca F)$ for any $\ca O_X$-module $\ca F$.\\

\subsection{The Hodge decomposition}
\
\indent The tools we are going to use are close to those of Swan \cite[2]{Swa}. To get a Hodge decompostion of higher order Hochschild cohomology of schemes, we also refer to Pirashvili's work \cite{Pir}.\\

For any ringed space $X$ and any pointed simplicial set $K$, one can define the \textit{$K$-homology sheaf of $X$} by $\ca H(K,X)=H\big(\ca C(K,X))$. We are particularly interested in the homology of the higher order Hochschild sheaf.\\

\noindent\hypertarget{4.1.1}{\textbf{Definition 4.1.1.}} Let $d$ be a positive integer and $X$ a ringed space. We define
\begin{center}
$\ca H^\D_X=H\big(\ca C_X^\D\big)$
\end{center}
the \textit{Higher Hochschild homology sheaf of order $d$ of $X$}.\\

We will use the notation $\ca H_X=\ca H^{[1]}_X$ for the Hochschild homology sheaf of $X$.\\

Just like Swan explained \cite[2.4]{Swa}, when $X$ is a scheme and $K$ is connected, the sheaves $\ca H_q(K,X)$ are quasi-coherent such that on each affine open set $U$ of $X$, we have
\begin{center}
$\ca H(K,X)(U)=H\big(K,\ca O_X (U)\big)$
\end{center}
This is a consequence of \hyperlink{3.2.2}{Theorem 3.2.2} and the fact that sheafification commutes with homology \cite[3.1]{Toh} so that the sheaf $\ca H(K,X)$ is associated to the presheaf
\begin{center}
$U\map H\big(K,\ca O_X (U)\big)$
\end{center}
One can compute the first higher order Hochschild homology sheaves.\\

\noindent\hypertarget{4.1.2}{\textbf{Example 4.1.2.}} Let $X$ be a separated scheme over a \textbf{field of any characteristic}. For any positive integer $d$, there is a natural isomorphism
\begin{center}
$\ca H^\D_{X,q}\iso{\left\{
    \begin{array}{ll}
        \ca O_X & \mbox{if } q=0\\
        0 & \mbox{if } 0<q<d\\
        \Omega^1_X & \mbox{if } q=d
    \end{array}
\right.}$
\end{center}
It was constructed on any affine open set $U$ of $X$ in \hyperlink{3.1.2}{Example 3.1.2}, since we always have
\begin{center}
$\Omega^1_X|_U\iso\Omega^1_U$
\end{center}
\cite[IV, 17.1.3.(i), 17.2.4]{EGAIV4}. These natural isomorphisms are compatible with restriction, which allows us to patch them \cite[II, Ex, 1.15]{Har}.\\

We can also compute higher degrees in the smooth case.\\

\noindent\hypertarget{4.1.3}{\textbf{Theorem 4.1.3.}} Let $d$ be a positive integer and $X$ a separated smooth scheme. There is a natural isomorphism
\begin{center}
$\ca H^\D_{X,q}\iso {\left\{
    \begin{array}{ll}
        \Omega^j_X & \mbox{if } q=jd\ \& \ d \mbox{ is odd}\\
        \ba{Sym}^j_{\ca O_X} \Omega^1_X & \mbox{if } q=jd\ \& \ d \mbox{ is even}\\
        0 & \mbox{if } d\nmid q
    \end{array}
\right.}$
\end{center}
\ \\
\textit{Proof.} We use the same argument as for \hyperlink{4.1.2}{Example 4.1.2} : this natural isomorphism is given on any affine open set $U$ of $X$ by \hyperlink{3.1.3}{Theorem 3.1.3} \cite[IV, 17.3.2.(ii)]{EGAIV4}. We just have to show that the restriction functor $|_U$ always commutes with exterior and symmetric powers : for any $\ca O_X$-module $\ca F$ and any positive integer $j$, there are colimit descriptions
\begin{center}
$\bigwedge^j_{\ca O_X}\ca F=\underset{\sigma\in\mathfrak{S}_j}{\ba{colim}}\big(\ca F^{\otimes j}\xrightarrow{\pm\sigma_*}\ca F^{\otimes j}\big)$\hspace*{.5cm}$\&$\hspace*{.5cm}$\ba{Sym}^j_{\ca O_X}\ca F=\underset{\sigma\in\mathfrak{S}_j}{\ba{colim}}\big(\ca F^{\otimes j}\overset{\sigma_*}{\arr}\ca F^{\otimes j}\big)$
\end{center}
where the sign $\pm$ is the signature of the permutation, and restriction functor $|_U$ commutes with colimits and tensor products \cite[0, 4.4.3.1, 4.3.3.1]{EGA}.\\

We study the properties of the higher order Hochschild homology sheaf because of the next proposition, which is inspired by the \textit{Hodge spectral sequence} of Swan \cite[2.3]{Swa}.\\

\noindent\hypertarget{4.1.4}{\textbf{Proposition 4.1.4.}} Let $X$ be a ringed space, $K$ a pointed simplicial set and $\ca F$ an $\ca O_X$-module. There is a natural spectral sequence
\begin{center}
$E_2^{pq}=\ba{Ext}^p_{\ca O_X}\big((\ca H_q(K,X),\ca F\big)\implies H^{p+q}(K,X,\ca F)$
\end{center}
\ \\
\textit{Proof.} In the same way that we computed the Hochschild cohomology of order $1$ of $X$ with coefficients in $\ca F$, we can choose an injective resolution of $\ca O_X$-modules $\ca F\arr \ca I$ to get a fibrant remplacement $\ca I$ of the complex of $\ca O_X$-modules $\ca F$. Now we have a double complex $E_0^{pq}=\ba{Hom}_{\ca O_X}\big(\ca (C_q(K,X),\ca I^p\big)$ with total cohomology $H(K,X,\ca F)$. If we filter it by columns \cite[5.6.1]{Wei} then we get the announced spectral sequence since the first page $E_1$ is the vertical cohomology which is preserved by the exact functor $\ba{Hom}_{\ca O_X}(\ \ \ ,\ca I^p)$ and the second page $E_2$ is the horizontal cohomology which computes the ext functors :
\begin{center}
$E_1^{pq}=H^q\Big(\ba{Hom}_{\ca O_X}\big(\ca C(K,X),\ca I^p\big)\Big)\iso \ba{Hom}_{\ca O_X}\big(\ca H_q(K,X),\ca I^p\big)$
\end{center}
\begin{center}
$E_2^{pq}\iso H^p\Big(\ba{Hom}_{\ca O_X}\big(\ca H_q(K,X),\ca I\big)\Big)=\ba{Ext}^p_{\ca O_X}\big(\ca H_q(K,X),\ca F\big)$
\end{center}
\ \\
\indent Using the spectral sequence of \hyperlink{4.1.4}{Proposition 4.1.4} for $K=\mathbf{S}^d$, we obtain the \textit{Hodge decomposition for higher order Hochschild cohomology of smooth varieties}.\\

\noindent\hypertarget{4.1.5}{\textbf{Theorem 4.1.5.}} Let $d$ be a positive integer, $X$ a separated smooth scheme and $\ca F$ an $\ca O_X$-module. If $d$ is odd then we have a natural isomorphism
\begin{center}
$H_\D^n(X,\ca F)\iso\bigoplus\limits_{p+jd=n}H^p(X,\ca T_X^j\otimes\ca F)$
\end{center}
where $\ca T_X^j=(\Omega^j_X)^\vee$. If $d$ is even then we have a natural isomorphism
\begin{center}
$H_\D^n(X,\ca F)\iso\bigoplus\limits_{p+jd=n}H^p(X,\ca S_X^j\otimes\ca F)$
\end{center}
where $\ca S_X^j=(\ba{Sym}^j_{\ca O_X} \Omega^1_X)^\vee$.\\

\noindent\textit{Proof.} First, we show that the spectral sequence of \hyperlink{4.1.4}{Proposition 4.1.4}
\begin{center}
$E_2^{pq}=\ba{Ext}^p_{\ca O_X}\big(\ca H^\D_{X,q},\ca F\big)\implies H_{[d]}^{p+q}(\ca O_X,\ca F)$
\end{center}
degenerates at the page $E_2$. For $d=1$, this was done by Swan \cite[2.6]{Swa}, but this is also a consequence of Loday's \textit{$\lambda$-decomposition} \cite[4.5.10]{Lod}. It induces a decomposition $\ca C^{[1]}_X=\bigoplus_{i\geq 0}\ca C^{(i)}$ such that $H_q(\ca C^{(i)})\neq 0$ if and only if $q=i$ \cite[4.5.12, 3.4.4]{Lod}. Therefore the page $E_2$ is a direct sum of single columns and the differentials are zero. For $d\geq 2$, we can use \hyperlink{4.1.3}{Theorem 4.1.3}. We have $\ca H^\D_{X,q}\neq 0$ if and only if $d$ divides $q$ so at the page $E_2$, each non-zero row is surrounded by zero rows :
\begin{center}
$\xymatrix{ & \vdots & \vdots & \vdots & \vdots & \vdots & \\
\cdots & 0 & 0 & 0 & 0 & 0 & \cdots \\
\cdots & E_2^{p+2,jd} & E_2^{p+1,jd} & E_2^{p,jd} \ar[rru] & E_2^{p-1,jd} & E_2^{p-2,jd} & \cdots \\
\cdots & 0 \ar[rru] & 0 & 0 & 0 & 0 & \cdots \\
 & \vdots & \vdots & \vdots & \vdots & \vdots & }$
\end{center}
Hence the differentials are zero and we have a natural isomorphism
\begin{center}
$H_{[d]}^n(\ca O_X,\ca F)\iso\bigoplus\limits_{p+q=n}\ba{Ext}^p_{\ca O_X}\big(\ca H^\D_{X,q},\ca F\big)=\bigoplus\limits_{p+jd=n}\ba{Ext}^p_{\ca O_X}\big(\ca H^\D_{X,jd},\ca F\big)$
\end{center}
Now we want to compute $\ba{Ext}^p_{\ca O_X}\big(\ca H^\D_{X,jd},\ca F\big)$. We can use the usual \textit{Grothendieck spectral sequence} \cite[2.4.1]{Toh} of $\ba{Hom}_{\ca O_X}\big(\ca H^\D_{X,jd},\ \ \ \big)=\Gamma\circ\ca Hom_{\ca O_X}\big(\ca H^\D_{X,jd},\ \ \ \big)$ :
\begin{center}
$F_2^{pq}=H^p\Big(X,\ca Ext^q_{\ca O_X}\big(\ca H^\D_{X,jd},\ca F\big)\Big)\implies \ba{Ext}^{p+q}_{\ca O_X}\big(\ca H^\D_{X,jd},\ca F\big)$
\end{center}
Since $X$ is smooth, $\ca H^\D_{X,jd}$ is locally free of finite rank \cite[IV, 17.2.3]{EGAIV4} and we have $\ca Ext^q_{\ca O_X}\big(\ca H^\D_{X,jd},\ca F\big)=0$ for $q>0$ \cite[4.2.3]{Toh}. Thus, the specral sequence $F$ induces a natural isomorphism
\begin{center}
$\ba{Ext}^p_{\ca O_X}\big(\ca H^\D_{X,jd},\ca F\big)\iso H^p\Big(X,\ca Hom_{\ca O_X}\big(\ca H^\D_{X,jd},\ca F\big)\Big)$
\end{center}
Moreover, the natural morphism $(\ca H^\D_{X,jd})^\vee\otimes\ca F\arr\ca Hom_{\ca O_X}\big(\ca H^\D_{X,jd},\ca F\big)$ is an isomorphism \cite[0, 5.4.2.1]{EGA}, which ends the proof.\\

\noindent\hypertarget{4.1.6}{\textbf{Example 4.1.6.}} Let $d$ be a positive integer, $X$ a separated smooth scheme and $\ca F$ an $\ca O_X$-module. One can visualize the higher order Hochschild cohomology groups :
\begin{center}
$\xymatrix @!0 @R=9mm @C=17mm { H_\D^0(X,\ca F) & \iso & \Gamma\ca F & & & & & & \\
\vdots & & & & & & & & \\
H_\D^n(X,\ca F) & \iso & H^n(X,\ca F) & & & & & & 0<n<d \\
\vdots & & & & & & & & \\
H_\D^d(X,\ca F) & \iso & H^d(X,\ca F) & \oplus & \Gamma(\ca T_X\otimes\ca F) & & & & \\
\vdots & & & & & & & \\
H_\D^n(X,\ca F) & \iso & H^n(X,\ca F) & \oplus & \ \ \ H^{n-d}(X,\ca T_X\otimes\ca F) & & & & d<n<2d \\
\vdots & & & & & & & & \\
H_\D^{2d}(X,\ca F) & \iso & H^{2d}(X,\ca F) & \oplus & H^d(X,\ca T_X\otimes\ca F) & \oplus & \Gamma(\ca D_X^2\otimes\ca F) & & \\
\vdots & & & & & & & & }$
\end{center}
where $\ca T_X=\ca T_X^1=\ca S_X^1=(\Omega^1_X)^\vee$ and $\ca D_X^j$ is $\ca T_X^j$ for $d$ odd and $\ca S_X^j$ for $d$ even. Assuming that $d>\ba{dim}(X)$, we have for any integer $n\geq 0$ a natural isomorphism
\begin{center}
$H_\D^n(X,\ca F)\iso H^p(X,\ca D_X^j\otimes\ca F)$
\end{center}
where $n=jd+p$ is the \textit{Euclidean division} of $n$ by $d$ with remainder $p$ \cite[III, 2.7]{Har}.\\

\subsection{Affine case}
\
\indent The definition of the $K$-cohomology of a scheme lead us to the following question. Does this definition be equivalent to the algebraic definition of the Higher Hochschild cohomology \cite[3.9]{Gin} ? In other words, do we have an isomorphism 
\begin{center}
$H(K,\ba{Spec}(A),\widetilde M)\iso H(K,A,M)$
\end{center}
for any simplicial set $K$ and any algebra $A$ ? \hyperlink{3.2.1.3}{Exemple 3.2.1.3} shows that generaly the complexe of sheaves $\ca C(K,\ba{Spec}(A))$ is not weak equivalent to the complexe of sheaves associated to the Higher Hochschild complex $C(K,A)$.\\

Using \hyperlink{3.2.2}{Theorem 3.2.2}, we will show that the $K$-cohomology of an affine scheme is isomorphic to the $K$-cohomology of the underlying algebra if the simplicial set $K$ is connected. We will also compare the spectral sequence of \hyperlink{4.1.4}{Proposition 4.1.4} and the spectral sequence of Pirashvili \cite[2.4]{Pir}.\\

Let $X=\ba{Spec}(A)$ be an affine scheme. Recall the functor
\begin{center}
$\sim\ : A\ba{-mod}\arr\ca O_X\ba{-mod}$
\end{center}
which associates to each $A$-module $M$ a quasi-coherent $\ca O_X$-module $\widetilde M$ \cite[I, 1.3]{EGA}. The following \textit{Quillen adjunction} is the main tool.\\

\noindent\hypertarget{4.2.1}{\textbf{Proposition 4.2.1.}} Let $X=\ba{Spec}(A)$ be an affine scheme. The pair
\begin{center}
$\sim\ \ :\xymatrix @!0 @R=5mm @C=25mm {
    \ba{Ch}(A) \ar@<3pt>[r] & \ba{Ch}(\ca O_X) \ar@<3pt>[l]}:\ \Gamma$
\end{center}
is a Quillen adjunction between injective model structures. Let $C$ be a complex of $A$-modules and $M$ a $A$-module. There is a natural weak equivalence of chain complexes
\begin{center}
$\mathbb{R}\ba{Hom}_{\ba{Ch}(\ca O_X)}\big(\widetilde C,\widetilde M\big)\we\mathbb{R}\ba{Hom}_{\ba{Ch}(A)}(C,M)$
\end{center}
\ \\
\noindent\textit{Proof.} It is a standard fact that this pair is an adjunction \cite[II, Ex, 5.3]{Har}. In order to show that it is a Quillen adjunction, we just have to show that the functor $\sim$ preserves cofibrations and trivial cofibrations \cite[1.3.4]{Hov1}, which is a consequence of its exactness \cite[II, 5.2]{Har}. Deriving the enriched adjunction $\ba{Hom}_{\ba{Ch}(\ca O_X)}\circ(\sim\times 1)\iso\ba{Hom}_{\ba{Ch}(\ca O_X)}\circ(1\times\Gamma)$ \cite[1.3.7]{Hov1}, we obtain the announced weak equivalence because functor $\sim$ is its left derived functor since all objects are cofibrant in $\ba{Ch}(A)$, and we have a natural weak equivalence of chain complexes of $A$-modules $M\we\mathbb{R}\Gamma\widetilde M$ since quasi-coherent $\ca O_X$-modules have no sheaf cohomology on $X=\ba{Spec}(A)$ \cite[III, 1.3.1]{EGA}.\\

Before getting to the main theorem of this section, we have to compare \hyperlink{4.0.1}{Definition 4.0.1} with a quasi-coherent version.\\

\newpage

\noindent\hypertarget{4.2.2}{\textbf{Lemma 4.2.2.}} Let $K$ a connected pointed simplicial set and $X=\ba{Spec}(A)$ an affine scheme. The canonical morphism of sheaves of CDGA on $X$
\begin{center}
$\widetilde{C(K,A)}\arr\ca C(K,X)$
\end{center}
is a quasi-isomorphism.\\

\noindent\textit{Proof.} This morphism is given on open subset $D(s)=\{\mathfrak{p}\in X : s\notin\mathfrak{p}\}$ of $X$ by
\begin{center}
$C(K,A)_s\arr C(K,A_s)\arr\ca C(K,X)\big(D(s)\big)$
\end{center}
\hyperlink{3.2.2}{Theorem 3.2.2} shows that the presheaves of CDGA on $X$
\begin{center}
$D(s)\map C(K,A)_s$\hspace*{5mm}$\&$\hspace*{5mm}$D(s)\map C(K,A_s)$
\end{center}
are quasi-isomorphic. Thus, the associated sheaves are quasi-isomorphic \cite[3.1]{Toh}.\\

\noindent\hypertarget{4.2.3}{\textbf{Theorem 4.2.3.}} Let $K$ be a connected pointed simplicial set, $X=\ba{Spec}(A)$ an affine scheme and $\ca F=\widetilde M$ a quasi-coherent $\ca O_X$-module. There is a natural isomorphism
\begin{center}
$H^n(K,X,\ca F)\iso H^n(K,A,M)$
\end{center}
\ \\
\noindent\textit{Proof.} \hyperlink{4.2.2}{Lemma 4.2.2} and \hyperlink{4.2.1}{Proposition 4.2.1} give us two natural weak equivalences
\begin{center}
$\mathbb{R}\ba{Hom}_{\ba{Ch}(\ca O_X)}\big(\ca C(K,X),\widetilde M\big)\we \mathbb{R}\ba{Hom}_{\ba{Ch}(\ca O_X)}\big(\widetilde{C(K,A)},\widetilde M\big)\we \mathbb{R}\ba{Hom}_{\ba{Ch}(A)}\big(C(K,A),M\big)$
\end{center}
and the derived hom complex $\mathbb{R}\ba{Hom}_{\ba{Ch}(A)}\big(C(K,A),M\big)$ can be computed by the chain complex $\ba{Hom}_{\ba{Ch}(A)}\big(C(K,A),M\big)$. In fact, for any integer $n\geq 0$, $A^{\otimes K_n\setminus\ast}$ is a vector space, so $A^{\otimes K_n}$ is a free $A$-module. This implies that $C(K,A)$ is a cofibrant complex of $A$-modules for the projective model structure \cite[2.3.6]{Hov1}.\\

\noindent\hypertarget{4.2.4}{\textbf{Remark 4.2.4.}} \textit{Projective schemes.} Let $X=\ba{Proj}(A)$ be a projective scheme associated to a graded algebra $A=\bigoplus_{n\in\mathbb N}A_n$ \cite[II, 2]{EGA}. As in \hyperlink{4.2.2}{Lemma 4.2.2}, for any connected pointed simplicial set $K$, the canonical morphism of sheaves of CDGA on $X$
\begin{center}
$\widetilde{C(K,A)}\arr\ca C(K,X)$
\end{center}
is a quasi-isomorphism, where $\widetilde{C_q(K,A)}$ is the $\ca O_X$-module associated to the graded $A$-module $C_q(K,A)=\bigoplus_{n\in\mathbb N}C_q(K,A_n)$ for any integer $q\geq 0$. In fact, this morphism is given on open subset $D_+(s)=\{\mathfrak{p}\in X : s\notin\mathfrak{p}\}$ of $X$ by
\begin{center}
$\big(C(K,A)_s\big)_0\arr C\big(K,(A_s)_0\big)\arr\ca C(K,X)\big(D_+(s)\big)$
\end{center}
so we just have to check that Hochschild differentials preserve the graduation of $A$, since they already commute with localization : for any map $f:L\arr L'$, we have
\begin{center}
$\Bigg|f_*\Big(\underset{x\in L}{\bigotimes}a_x\Big)\Bigg|=\Bigg|\underset{y\in L'}{\bigotimes}\Big(\underset{f(x)=y}{\prod}a_x\Big)\Bigg|=\underset{y\in L'}{\sum}\Big(\underset{f(x)=y}{\sum}|a_x|\Big)=\underset{x\in L}{\sum}|a_x|=\Bigg|\underset{x\in L}{\bigotimes}a_x\Bigg|$
\end{center}
where $|a|$ is the degree of an homogeneous element $a$ of any graded algebra. As in \hyperlink{4.2.1}{Proposition 4.2.1}, we also have a Quillen adjonction
\begin{center}
$\sim\ \ :\xymatrix @!0 @R=5mm @C=35mm {
    \ba{Ch}(A-\ba{grmod}) \ar@<3pt>[r] & \ba{Ch}(\ca O_X) \ar@<3pt>[l]}:\ \Gamma_*$
\end{center}
but in general, we don't have an analogue of the weak equivalence between the derived hom complexes as in \hyperlink{4.2.1}{Proposition 4.2.1}, since quasi-coherent $\ca O_X$-modules may have sheaf cohomology on $X$ \cite[III, 2.2.2]{EGA}.\\

\newpage

Following \hyperlink{4.2.3}{Theorem 4.2.3}, we want to compare the spectral sequences approaching $K$-cohomology of affine schemes. Let $\ba{Vect}_\ba{Fin'}$ be the category of contravariant functors over the category of pointed finite sets with values in the category of vector spaces. Pirashvili associated to any functor $F:\big(\ba{Fin}'\big)^\ba{op}\arr\ba{Vect}$ and any pointed finite simplicial set $K$ a cosimplicial vector space $F(K)$ by taking the composition
\begin{center}
$F(K):\xymatrix{\mathbf{\Delta} \ar[r]^-K & \big(\ba{Fin}'\big)^\ba{op} \ar[r]^-F & \ba{Vect}}$
\end{center}
and then a spectral sequence approaching its cohomology
\begin{center}
$E_2^{pq}=\ba{Ext}^p_{\ba{Vect}_\ba{Fin'}}\Big(\ca J_q\big(H(K)\big),F\Big)\implies H^{p+q}\big(F(K)\big)$
\end{center}
where $\ca J$ is the graded version of the dual Loday functor \cite[1.7]{Pir} and $H(K)$ is the homology of $K$ with coefficients in $k$. We are particularly interested in the case where $K$ is connected and $F$ is the functor $\ca H(A,M)=\ba{Hom}_A\big(\ca L(A,A),M\big)$ for an algebra $A$ and an $A$-module $M$ so that the cohomology of $F(K)$ is $K$-cohomology of $A$ with coefficients in $M$ \cite[3.9]{Gin}.\\

Let $\ba{Ch}_\ba{Fin'}$ be the category of contravariant functors over the category of pointed finite sets with values in the category of chain complexes. It is a model category such that the weak equivalences are the quasi-isomorphisms and the fibrations are the epimorphisms \cite[11.6.1]{Hir}. The cofibrant objects of $\ba{Ch}_\ba{Fin'}$ are \textit{projective in the sense of Spaltenstein} \cite[1.4]{Spa} \cite[1.2.10]{Hov1} so the hom complex functor
\begin{center}
$\ba{Hom}_{\ba{Ch}_\ba{Fin'}}:\big(\ba{Ch}_\ba{Fin'}\big)^\ba{op}\times\ba{Ch}_\ba{Fin'}\arr\ba{Ch}$
\end{center}
has a right derived functor. To get an isomorphism between the spectral sequence of \hyperlink{4.1.4}{Proposition 4.1.4} and the Pirashvili's spectral sequence, we will use another Quillen adjunction.\\

\noindent\hypertarget{4.2.5}{\textbf{Proposition 4.2.5.}} Let $A$ be an algebra. The pair
\begin{center}
$\otimes_{\ba{Fin}'}\ca L(A,A)\ :\xymatrix @!0 @R=5mm @C=25mm {
    \ba{Ch}_{\ba{Fin}'} \ar@<3pt>[r] & \ba{Ch}(A) \ar@<3pt>[l]}:\ \ca H(A,\ \ )$
\end{center}
is a Quillen adjunction. Let $F:\big(\ba{Fin}'\big)^\ba{op}\arr\ba{Vect}$ be a functor and $M$ a $A$-module. There is a natural weak equivalence of chain complexes
\begin{center}
$\mathbb{R}\ba{Hom}_{\ba{Ch}(A)}\Big(F\overset{\mathbb{L}}{\underset{\ba{Fin}'}{\otimes}}\ca L(A,A),M\Big)\we\mathbb{R}\ba{Hom}_{\ba{Ch}_{\ba{Fin}'}}\big(F,\ca H(A,M)\big)$
\end{center}
\ \\
\textit{Proof.} It is a standard fact that this pair is an adjunction \cite[IX, 5, 6]{ML2}. In order to show that it is a Quillen adjunction, we just have to show that functor $\ca H(A,\ \ )$ preserves fibrations and trivial fibrations \cite[1.3.4]{Hov1}, which is a consequence of its exactness. In fact, we saw that for any positive integer $n$, $A^{\otimes n-1}$ is a vector space, so $A^{\otimes n}$ is a free $A$-module. Deriving the enriched adjunction $\ba{Hom}_{\ba{Ch}(A)}\circ\big(\otimes_{\ba{Fin}'}\ca L(A,A)\times 1\big)\iso\ba{Hom}_{\ba{Ch}_{\ba{Fin}'}}\circ\big(1\times\ca H(A,\ \ )\big)$ \cite[1.3.7]{Hov1}, we obtain the announced weak equivalence because $\ca H(A,\ \ )=\mathbb{R}\ca H(A,\ \ )$ since all object of $\ba{Ch}(A)$ are fibrant.\\

Remark that we could have formulated a more general adjunction with an arbitrary functor $G:\ba{Fin}'\arr\ba{Vect}$ instead of $\ca L(A,A)$, so that derived functor $\otimes^\mathbb{L}_\ba{Fin'}G$ always exists. For any pointed finite simplicial set $L$, Pirashvili defined a functor
\begin{center}
$h_L:(\ba{Fin}')^\ba{op}\arr\ba{sVect}$
\end{center}
which associates to each pointed finite set $K$ the \textit{free vector space} over $\ba{Hom}_{\ba{Fin}'}(K,L)$. The simplicial structure of $h_L(K)$ is naturally given by the simplicial structure of $L$. Pirashvili showed that we have a natural isomorphism of simplicial vector spaces
\begin{center}
$h_L\otimes_{\ba{Fin}'}G\iso G(L)$
\end{center}
\cite[1.5]{Pir}. We can compute the homology of $G(L)$ with this formula. \\

\noindent\hypertarget{4.2.6}{\textbf{Lemma 4.2.6.}} Let $q$ be a non-negative integer, $K$ a pointed finite simplicial set and $F:\ba{Fin'}\arr\ba{Vect}$ a functor. There is a natural weak equivalence of chain complexes
\begin{center}
$H_q(h_K)\overset{\mathbb{L}}{\underset{\ba{Fin}'}{\otimes}}F\we H_q\big(F(K)\big)$
\end{center}
\ \\
\textit{Proof.} $H_q(h_K)\overset{\mathbb{L}}{\underset{\ba{Fin}'}{\otimes}}F\we H_q\Big(h_K\overset{\mathbb{L}}{\underset{\ba{Fin}'}{\otimes}}F\Big)=H_q\big(h_K\otimes_{\ba{Fin}'}F\big)\iso H_q\big(F(K)\big)$ \cite[1.1]{Pir}.\\

\noindent\hypertarget{4.2.7}{\textbf{Theorem 4.2.7.}} Let $K$ be a connected pointed simplicial set, $X=\ba{Spec}(A)$ an affine scheme and $\ca F=\widetilde M$ a quasi-coherent $\ca O_X$-module. The following spectral sequences are naturally isomorphic
\begin{center}
$\ba{Ext}^p_{\ca O_X}\big(\ca H_q(K,X),\ca F\big)\implies H^{p+q}(K,X,\ca F)$
\end{center}
\begin{center}
$\ba{Ext}^p_{\ba{Vect}_\ba{Fin'}}\Big(\ca J_q\big(H(K)\big),\ca H(A,M)\Big)\implies H^{p+q}(K,A,M)$
\end{center}
\ \\
\textit{Proof.} \hyperlink{4.2.2}{Lemma 4.2.2} and \hyperlink{4.2.1}{Proposition 4.2.1} give rise to natural isomorphisms
\begin{center}
$\ba{Ext}^p_{\ca O_X}\big(\ca H_q(K,X),\ca F\big)\iso \ba{Ext}^p_{\ca O_X}\Big(\widetilde{H_q(K,A)},\widetilde M\Big)\iso \ba{Ext}^p_A\big(H_q(K,A),M\big)$
\end{center}
for any integers $p,q\geq 0$. \hyperlink{4.2.6}{Lemma 4.2.6} and \hyperlink{4.2.5}{Proposition 4.2.5} give us natural isomorphisms
\begin{center}
$\ba{Ext}^p_A\big(H_q(K,A),M\big)=H^p\Big(\mathbb{R}\ba{Hom}_{\ba{Ch}(A)}\big(H_q(K,A),M\big)\Big)$
\end{center}
\begin{center}
$\iso H^p\Bigg(\mathbb{R}\ba{Hom}_{\ba{Ch}(A)}\Big(H_q(h_K)\overset{\mathbb{L}}{\underset{\ba{Fin}'}{\otimes}}\ca L(A,A),M\Big)\Bigg)$
\end{center}
\begin{center}
$\iso H^p\Big(\mathbb{R}\ba{Hom}_{\ba{Ch}_{\ba{Fin}'}}\big(H_q(h_K),\ca H(A,M)\big)\Big)=\ba{Ext}^p_{\ba{Vect}_\ba{Fin'}}\big(H_q(h_K),\ca H(A,M)\big)$
\end{center}
for any integer $p,q\geq 0$, which ends the proof since Pirashvili showed that the second spectral sequence is in fact isomorphic to the spectral sequence
\begin{center}
$\ba{Ext}^p_{\ba{Vect}_\ba{Fin'}}\big(H_q(h_K),\ca H(A,M)\big)\implies H^{p+q}\Big(\ba{Hom}_{\ba{Vect}_\ba{Fin'}}\big(h_K,\ca H(A,M)\big)\Big)$
\end{center}
\cite[1.6, 2.4]{Pir}.\\

\noindent\hypertarget{4.2.8}{\textbf{Corollary 4.2.8.}} Let $d$ be a positive integer, $X=\ba{Spec}(A)$ an affine scheme and $\ca F=\widetilde M$ a quasi-coherent $\ca O_X$-module. The following spectral sequences are naturally isomorphic
\begin{center}
$\ba{Ext}^p_{\ca O_X}\big(\ca H^\D_{X,q},\ca F\big)\implies H^{p+q}_\D (X,\ca F)$
\end{center}
\begin{center}
$\ba{Ext}^p_{\ba{Vect}_\ba{Fin'}}\Big(\ca J_q\big(H(\mathbf{S}^d)\big),\ca H(A,M)\Big)\implies H^{p+q}_\D (A,M)$
\end{center}

\hypertarget{5}{\section{Gerstenhaber-Schack type definition}}
\
\indent In all of this section, if nothing else is specified, $k$ is a field of caracteristic zero. Recall that algebras are associatives and commutatives over $k$ with identity.\\

Let $X$ be a separated scheme. We call the contravariant functors over the category $\ba{Aff}(X)$ of affine open sets of $X$ \textit{presheaves on $X$}. One can also build a Higher Hochschild functor with presheaves by taking for any finite simplicial set $K$ and any presheaf of algebras $\ca O$ on $X$ the presheaf of CDGA on $X$
\begin{center}
$C(K,\ca O):U\arr C\big(K,\ca O(U)\big)$
\end{center}
We can define $C^\D(\ca O)=C(\mathbf{S}^d,\ca O)$ and $H^\D(\ca O)=H\big(C^\D(\ca O)\big)$ as in \hyperlink{3.1.1}{Definition 3.1.1} for any positive integer $d$. \hyperlink{3.1.2}{Example 3.1.2} shows that we always have an isomorphism of algebras $H^\D_0(\ca O)\iso\ca O$ which allows us to think of any $\ca O$-module as a $C^\D(\ca O)$-module via the morphism of presheaves of CDGA on $X$
\begin{center}
$C^\D(\ca O)\arr\ca O$
\end{center}
Following the \textit{Gerstenhaber-Schack} type definition \cite[3]{Swa} and the weak equivalence of simplicial set $\mathbf{S}^d\we\ast\cup^h_{\mathbf{S}^{d-1}}\ast$, we would like to define the \textit{Higher Hochschild cohomology of order $d$ of $X$ with coefficients in any $C^{[d-1]}(\ca O)$-module $\ca M$} by taking the cohomology of $\mathbb{R}\ba{Hom}_{C^{[d-1]}(\ca O)}(\ca O,\ca M)$ where $\ca O$ is the structure presheaf of $X$ and with the convention $C^{[0]}(\ca O)=\ca O^e=\ca O\otimes\ca O$. We therefore need to give a meaning to this derived hom complex.\\

\hypertarget{5.1}{\subsection{Model structure of modules over a presheaf of CDGA}}
\
\indent Throughout this paragraph, $k$ may be a field of \textbf{any caracteristic}.\\

The category $\ba{Ch}(X)$ of presheaves of chain complexes over a separated scheme $X$ is a \textit{closed symmetric monoidal category} with the pointwise tensor product $\otimes$ and the hom presheaf $\ca Hom$. It is also a cofibrantly generated model category such that the generating cofibrations are the maps $(S^{n-1})_U\arr (D^n)_U$ and the generating trivial cofibrations are the maps $0\arr (D^n)_U$ \cite[11.6.1]{Hir}. Here we follow notations of Hovey \cite[2.3.3]{Hov1} and we denote by $C_U$ the constant presheaf on $U$ with value $C$ extended by zero on $X$ \cite[II, Ex, 1.19]{Har}. Remark that there is a natural isomorphism $\ba{Hom}_{\ba{Ch}(X)}\big((S^n)_U,\ca D\big)\iso\ca D_n(U)$ so that $(S^n)_U$ is a \textit{small object} of $\ba{Ch}(X)$ \cite[2.1.3]{Hov1}.\\

The following construction can be done with any presheaf of DGA on a basis of open sets of a topological space. Because of our context, we prefere to work with a presheaf of CDGA on the affine open sets of a separated scheme.\\

\noindent\hypertarget{5.1.1}{\textbf{Theorem 5.1.1.}} Let $\ca A$ be a presheaf of CDGA on a separated scheme $X$. The category $\ca A\ba{-mod}$ of presheaves of $\ca A$-modules on $X$ is a cofibrantly generated monoidal model category such that the weak equivalences are the quasi-isomorphisms and the fibrations are the epimorphisms. Moreover, if $\ca M\arr\ca N$ is a cofibration of $\ca A$-modules then $\ca M(U)\arr\ca N(U)$ is a cofibration of $\ca A(U)$-modules for any affine open set $U$ of $X$.\\

\noindent\textit{Proof.} The category $\ca A\ba{-mod}$ can be thought of as the category of algebras over the monad
\begin{center}
$\ca A\otimes :\ba{Ch}(X)\arr\ba{Ch}(X)$
\end{center}
\cite[VI]{ML2}. This functor commutes with colimits, and any object of $\ba{Ch}(X)$ is fibrant, so we just have to show that every $\ca A$-module has a \textit{path object} to obtain the model structure \cite[3.1, A.3]{SS}. Let $I$ be the bounded below chain complex
\begin{center}
$\xymatrix @!0 @R=5mm @C=18mm {\cdots \ar[r] & 0 \ar[r] & k \ar[r] & k\oplus k\\
 & & 1 \ar@{|->}[r] & (1,-1)}$
\end{center}
For any chain complex $C$, the hom complex $C^I=\ba{Hom}(I,C)$ is a path object for $C$. Since weak equivalences and fibrations are defined pointwise in $\ba{Ch}(X)$, the presheaf $\ca M^I:U\map \ca M(U)^I$ is a path object for any $\ca A$-module $\ca M$. For cofibrations, consider the product of section functors over affine open sets of $X$
\begin{center}
$\Gamma:\ca A\ba{-mod}\arr\underset{U\in\ba{Aff}(X)}{\prod}\big(\ca A(U)\ba{-mod}\big)$
\end{center}
It has a right adjoint which associates to each family $(M_U)_{U\in\ba{Aff}(X)}$ the $\ca A$-module
\begin{center}
$U\map\underset{V\subset U}{\prod}M_V$
\end{center}
This right adjoint functor is exact, so $\Gamma$ is in fact a left Quillen functor \cite[1.3.4]{Hov1}.\\

Now we can define the derived tensor product and the derived hom complex over a presheaf of CDGA on a separated scheme \cite[4.3.1]{Hov1}. The derived tensor product can be computed on any affine open set in the same way as the classical tensor product.\\

\noindent\hypertarget{5.1.2}{\textbf{Lemma 5.1.2.}} Let $\ca A$ be a presheaf of CDGA on a separated scheme $X$, $\ca M$ and $\ca N$ two $\ca A$-modules. For any affine open set $U$ of $X$, we have
\begin{center}
$\Big(\ca M\overset{\mathbb{L}}{\underset{\ca A}{\otimes}}\ca N\Big)(U)=\ca M(U)\overset{\mathbb{L}}{\underset{\ca A(U)}{\otimes}}\ca N(U)$
\end{center}
\ \\
\textit{Proof.} By \hyperlink{5.1.1}{Theorem 5.1.1}, if $\ca Q$ is a cofibrant remplacement of the $\ca A$-module $\ca M$ then $\ca Q(U)$ is a cofibrant remplacement of the $\ca A(U)$-module $\ca M(U)$ and we have
\begin{center}
$\big(\ca Q\otimes_{\ca A}\ca N\big)(U)=\ca Q(U)\otimes_{\ca A(U)}\ca N(U)$
\end{center}
\ \\

\indent We conclude this section by giving a derived base change lemma.\\

\noindent\hypertarget{5.1.3}{\textbf{Lemma 5.1.3.}} Let $\ca A\arr\ca B$ be a morphism of presheaves of CDGA on a separated scheme $X$, $\ca M$ a $\ca A$-module and $\ca N$ a $\ca B$-module. There is a natural weak equivalence of chain complexes
\begin{center}
$\mathbb{R}\ba{Hom}_{\ca B}\Big(\ca B\overset{\mathbb{L}}{\underset{\ca A}{\otimes}}\ca M,\ca N\Big)\we\mathbb{R}\ba{Hom}_{\ca A}(\ca M,\ca N)$
\end{center}
\ \\
\textit{Proof.} In fact, one can obtain a weak equivalence of presheaves of chain complexes on $X$ deriving the enriched adjunction $\ca Hom_{\ca B}\circ(\ca B\otimes_{\ca A}\times 1)\iso\ca Hom_{\ca A}$ \cite[IX, 5, 6]{ML2} \cite[1.3.7]{Hov1}, because the forgetful functor is its right dervied functor since all objects are fibrant in $\ca A$-mod.\\

\noindent\hypertarget{5.1.4}{\textbf{Remark 5.1.4.}} Let $\ca A$ be a presheaf of CDGA on a separated scheme $X$. One can see that if $\ca M$ is a cofibrant $\ca A$-module then $\ca M_x$ is a cofibrant $\ca A_x$-module for any $x\in X$. In fact, if $f:Y\arr X$ is an affine map and $\ca A\arr f_*\ca B$ a morphism of presheaves of CDGA on $X$, then the pair
\begin{center}
$f^*\ :\xymatrix @!0 @R=5mm @C=25mm {
    \ca A\ba{-mod} \ar@<3pt>[r] & \ca B\ba{-mod} \ar@<3pt>[l]}:\ f_*$
\end{center}
is a Quillen adjunction, since $f_*$ is always an exact functor \cite[1.3.4]{Hov1}. We can apply this to $x:\ast\arr X$ and $\ca A\arr x_*\ca A_x$.\\

\subsection{Link between both definitions}
\
\indent \hyperlink{5.1}{Section 5.1} gives rise to the following definition.\\

\noindent\hypertarget{5.2.1}{\textbf{Definiton 5.2.1.}} Let $d$ be a positive integer, $X$ a separated scheme with structure presheaf $\ca O$ and $\ca F$ a $\ca O$-module. We define
\begin{center}
$H_\D (\ca O,\ca F)=H\big(\mathbb{R}\ba{Hom}_{C^{[d-1]}(\ca O)}(\ca O,\ca F)\big)$
\end{center}
the \textit{Higher Hochschild cohomology of order $d$ of $\ca O$ with coefficients in $\ca F$}.\\

We recover the \textit{Gerstenhaber-Schack} type definiton for $d=1$ \cite[3]{Swa}. Recall our convention $C^{[0]}(\ca O)=\ca O^e$ and choose a projective resolution of $\ca O^e$-module $\ca P\arr\ca F$ \cite[1.10.1]{Toh}. Since fibrations are epimorphisms in $\ba{Ch}(\ca O^e)$, $\ca P$ is a cofibrant remplacement of the complex of $\ca O^e$-module $\ca F$. Hence the chain complex $\ba{Hom}_{\ba{Ch}(\ca O^e)}(\ca P,\ca F)$ computes both the derived hom complex $\mathbb{R}\ba{Hom}_{\ba{Ch}(\ca O^e)}(\ca O,\ca F)$ and the ext functor $\ba{Ext}_{\ca O^e}(\ca O,\ca F)$.\\

To get a link between \hyperlink{4.0.1}{Definition 4.0.1} and \hyperlink{5.2.1}{Definition 5.2.1} for a separated scheme, we need to pass from a presheaf to a sheaf and vice versa. Denote by $+$ the sheafification functor and by $\#$ the forgetful functor.\\

\noindent\hypertarget{5.2.2}{\textbf{Lemma 5.2.2.}} Let $X$ be a separated scheme with structure presheaf $\ca O$. The pair
\begin{center}
$+\ :\xymatrix @!0 @R=5mm @C=25mm {
    \ba{Ch}(\ca O) \ar@<3pt>[r] & \ba{Ch}(\ca O_X) \ar@<3pt>[l]}:\ \#$
\end{center}
is a Quillen adjunction.\\

\noindent\textit{Proof.} It is a standard fact that this pair is an adjunction \cite[II, 1.2]{Har}. In order to show that it is a Quillen adjunction, we just have to show that functor $+$ preserves cofibrations and trivial cofibrations \cite[1.3.4]{Hov1}, which is a consequence of \hyperlink{5.1.4}{Remark 5.1.4} \cite[2.3.9]{Hov1} and its exactness \cite[3.1]{Toh}.\\

Now we can prove that \hyperlink{4.0.1}{Definition 4.0.1} and \hyperlink{5.2.1}{Definition 5.2.1} coincide for separated schemes and quasi-coherent sheaves.\\

\noindent\hypertarget{5.2.3}{\textbf{Theorem 5.2.3.}} Let $d$ be a positive integer, $X$ a separated scheme with structure presheaf $\ca O$ and $\ca F$ an $\ca O_X$-module with no sheaf cohomology on affine open sets of $X$. There is a natural isomorphism
\begin{center}
$H^n_\D(X,\ca F)\iso H^n_\D(\ca O,\ca F^\#)$
\end{center}
\ \\
\textit{Proof.} Deriving the enriched adjunction $\ba{Hom}_{\ba{Ch}(\ca O_X)}\circ(+\times 1)\iso\ba{Hom}_{\ba{Ch}(\ca O)}\circ(1\times\#)$ \cite[1.3.7]{Hov1}, we have a natural weak equivalence of chain complexes
\begin{center}
$\mathbb{R}\ba{Hom}_{\ba{Ch}(\ca O_X)}\big(\ca C^\D_X,\ca F\big)\we\mathbb{R}\ba{Hom}_{\ba{Ch}(\ca O)}\big(C^\D(\ca O),\ca F^\#\big)$
\end{center}
In fact, functors $+$ and $\#$ can be derived by \hyperlink{5.2.2}{Lemma 5.2.2} and we have a natural weak equivalence $\mathbb{L}+\we +$ by exactness \cite[3.1]{Toh}. Moreover, $\ca C^\D_X=C^\D(\ca O)^+$ and we have a natural weak equivalence of complexes of $\ca O_X$-modules $\ca F\we \ca F^{\mathbb{R}\#}$ by hypothesis. \hyperlink{5.1.2}{Lemma 5.1.2} and the weak equivalence of simplicial set $\mathbf{S}^d\we\ast\cup^h_{\mathbf{S}^{d-1}}\ast$ give rise to a natural weak equivalance of $\ca O$-modules
\begin{center}
$C^\D(\ca O)\we\ca O\overset{\mathbb{L}}{\underset{C^{[d-1]}(\ca O)}{\otimes}} \ca O$
\end{center}

\newpage

\noindent and so we have natural weak equivalences of chain complexes
\begin{center}
$\mathbb{R}\ba{Hom}_{\ba{Ch}(\ca O)}\big(C^\D(\ca O),\ca F^\#\big)\we\mathbb{R}\ba{Hom}_{\ba{Ch}(\ca O)}\Big(\ca O\overset{\mathbb{L}}{\underset{C^{[d-1]}(\ca O)}{\otimes}} \ca O,\ca F^\#\Big)\we\mathbb{R}\ba{Hom}_{C^{[d-1]}(\ca O)}(\ca O,\ca F^\#)$
\end{center}
by \hyperlink{5.1.3}{Lemma 5.1.3}.\\

\noindent\hypertarget{5.2.4}{\textbf{Remark 5.2.4.}} Let $X$ be a scheme and $\ca F$ a quasi-coherent $\ca O_X$-module. We saw in the proof of \hyperlink{4.2.1}{Proposition 4.2.1} that $\ca F$ has no sheaf cohomology on affine open sets of $X$ \cite[II, 5.4]{Har} \cite[III, 1.3.1]{EGA}. Thus we can apply \hyperlink{5.2.3}{Theorem 5.2.3} when $X$ is separated.\\

\hypertarget{6}{\section{Swan definition of classical Hochschild cohomology}}
\
\indent In all of this section, $k$ is a field of \textbf{any characteristic}. Recall that algebras are associatives and commutatives over $k$ with identity.\\

Swan defined the Hochschild cohomology of a scheme $X$ with coefficient in an $\ca O_X$-module $\ca F$ by $\ba{Ext}_{\ca O_{X\times X}}(\delta_*\ca O_X,\delta_*\ca F)$ where $\delta :X\arr X\times X$ is the diagonal morphism of $X$. He showed that this ext functor coincides with the hyperext functor $\mathbb{E}\ba{xt}_{\ca O_X}(\ca C_X,\ca F)$ if $X$ is a \textit{quasi-projective} scheme \cite[2.1]{Swa}. We will generalize this fact and show that it is true for any separated scheme using \textit{derived tensor product of sheaves} and \textit{derived inverse/direct image} \cite[6, C, D]{Spa}.\\

Swan's definition comes with the usual \textit{Grothendieck spectral sequence} \cite[2.4.1]{Toh}
\begin{center}
$H^p\big(X\times X,\ca Ext^q_{\ca O_{X\times X}}(\delta_*\ca O_X,\delta_*\ca F)\big)\implies\ba{Ext}^{p+q}_{\ca O_{X\times X}}(\delta_*\ca O_X,\delta_*\ca F)$
\end{center}
We will show that it is isomorphic to the spectral sequence of \hyperlink{4.1.4}{Proposition 4.1.4} when, for example, $X$ is a separated smooth scheme.\\

\subsection{Strong Swan theorem}
\
\indent Let $X$ be a ringed space. The tensor product of complexes of $\ca O_X$-modules can be derived by means of a \textit{left Spaltenstein flat resolution} of either of the factors \cite[5.1, 5.9, 6.5]{Spa}. There is a relation between derived tensor products of sheaves and presheaves.\\

\noindent\hypertarget{6.1.1}{\textbf{Lemma 6.1.1.}} Let $\ca O$ be a presheaf of algebras on a separated scheme $X$, $\ca C$ and $\ca D$ two complexes of $\ca O^+$-modules. There is a natural weak equivalence of complexes of $\ca O^+$-modules
\begin{center}
$\Big(\ca C^\#\overset{\mathbb{L}}{\underset{\ca O}{\otimes}}\ca D^\#\Big)^+\we\ca C\overset{\mathbb{L}}{\underset{\ca O^+}{\otimes}}\ca D$
\end{center}
\ \\
\textit{Proof.} Let $\ca Q$ be a cofibrant remplacement of the complex of $\ca O$-modules $\ca C^\#$. The complex of $\ca O$-modules $\ca Q\otimes_{\ca O}\ca D^\#$ computes the derived tensor product $\ca C^\#\otimes^\mathbb{L}_{\ca O}\ca D^\#$. By \hyperlink{5.1.4}{Remark 5.1.4}, $\ca Q_x$ is cofibrant in $\ba{Ch}(\ca O_x)$ for any $x\in X$ so $\ca Q^+\arr\ca C^{\#+}\iso\ca C$ is a left Spaltenstein flat resolution over $\ca O^+$ \cite[5.5]{Hov2} \cite[5.3]{Spa} and the complex of $\ca O^+$-modules $\ca Q^+\otimes_{\ca O^+}\ca D$ computes the derived tensor product $\ca C\otimes^\mathbb{L}_{\ca O^+}\ca D$. Finally, we have an isomorphism of complexes of $\ca O^+$-modules
\begin{center}
$(\ca Q\otimes_{\ca O}\ca D^\#)^+\arr(\ca Q^{+\#}\otimes_{\ca O}\ca D^\#)^+\arr(\ca Q^{+\#}\otimes_{\ca O^{+\#}}\ca D^\#)^+=\ca Q^+\otimes_{\ca O^+}\ca D$
\end{center}
since $(\ca Q\otimes_{\ca O}\ca D^\#)^+_x\iso\ca Q^+_x\otimes_{\ca O^+_x}\ca D_x$ for any $x\in X$.\\

\newpage

\ \\
\indent Now let $f:X\arr Y$ be a morphism of ringed spaces. The functors $f^*$ and $f_*$ can also be derived by means of a certain type of resolution such that $\mathbb{L}f^*$ and $\mathbb{R}f_*$ come with derived enriched adjunctions \cite[6.7]{Spa}. Consider the case of an inclusion $i:Z\arr X$ of a closed set $Z$ of $X$. For any $\ca O_Z$-module $\ca F$ and any $x\in X$, there is a natural isomorphism of algebras
\begin{center}
$(i_*\ca F)_x\iso {\left\{
    \begin{array}{ll}
        \ca F_x & \mbox{if } x\in Z\\
        0 & \mbox{if } x\notin Z
    \end{array}
\right.}$
\end{center}
\cite[II, Ex, 1.19]{Har} and so the counit $i^{-1}i_*\arr 1$ is a natural isomorphism.\\

\noindent\hypertarget{6.1.2}{\textbf{Lemma 6.1.2.}} Let $i:X\arr Y$ be a morphism of ringed spaces which is an homeomorphism onto a closed set and $\ca C$ a complex of $\ca O_X$-modules. There is a natural weak equivalence of complexes of $\ca O_X$-modules
\begin{center}
$\mathbb{L}i^*i_*\ca C\we\ca O_X\overset{\mathbb{L}}{\underset{i^{-1}\ca O_Y}{\otimes}}\ca C$
\end{center}
\ \\
\textit{Proof.} Let $\ca Q\arr i_*\ca C$ be a left Spaltenstein flat resolution over $\ca O_Y$. The complex of $\ca O_X$-modules $i^*\ca Q$ computes the derived functor $\mathbb{L}i^*i_*\ca C$. Moreover, $i^{-1}\ca Q\arr i^{-1}i_*\ca C\iso\ca C$ is a left Spaltenstein flat resolution over $i^{-1}\ca O_Y$ \cite[5.3]{Spa} and so the complex of $\ca O_X$-modules $\ca O_X\otimes_{i^{-1}\ca O_Y}(i^{-1}\ca Q)$ computes the derived tensor product $\ca O_X\otimes^\mathbb{L}_{i^{-1}\ca O_Y}\ca C$. We can conclude the proof by definition of the inverse image functor $i^*\ca Q=\ca O_X\otimes_{i^{-1}\ca O_Y}(i^{-1}\ca Q)$.\\

\noindent\hypertarget{6.1.3}{\textbf{Proposition 6.1.3.}} Let $X$ be a scheme with diagonal morphism $\delta:X\arr X\times X$. There is a natural isomorphism of sheaves of algebras on $X$
\begin{center}
$\delta^{-1}\ca O_{X\times X}\iso{\ca O_X}^e$
\end{center}
\ \\
\textit{Proof.} The morphism is induced for any open sets $U$ and $V$ of X by the restriction maps
\begin{center}
$\ca O_X(U)\otimes\ca O_X(V)\arr\ca O_X(U\cap V)\otimes\ca O_X(U\cap V)$
\end{center}
Let $U=\ba{Spec}(A)$ be an affine open set of $X$. We have to identify the following two presheaves on $U$
\begin{center}
$D(r)\map\underset{D(s\otimes t)\supset\delta D(r)}{\ba{colim}}(A_s\otimes A_t)$\hspace*{1cm}$\&$\hspace*{1cm}$D(r)\map{A_r}^e$
\end{center}
where $D(r)=\{\mathfrak{p}\in U : r\notin\mathfrak{p}\}$ for any $r\in A$. Remark that
\begin{center}
$D(s\otimes t)\supset\delta D(r)\iff\sqrt{st}\ni r\implies \sqrt{s\otimes t}\ni r\otimes r\implies D(s\otimes t)\supset D(r\otimes r)\supset\delta D(r)$
\end{center}
so the functor $\underset{D(s\otimes t)\supset\delta D(r)}{\ba{colim}}$ is given by taking the value on $D(r\otimes r)$ \cite[IX, 3, Ex, 1]{ML2}.\\

We come now to the main theorem of this section. It generalizes a theorem of Swan which was initialy about quasi-projectif schemes \cite[2.1]{Swa}.\\

\noindent\hypertarget{6.1.4}{\textbf{Theorem 6.1.4.}} Let $X$ be a separated scheme with diagonal morphism $\delta:X\arr X\times X$ and $\ca F$ an $\ca O_X$-module. There is a natural isomorphism
\begin{center}
$HH(X,\ca F)\iso\ba{Ext}_{\ca O_{X\times X}}\big(\delta_*\ca O_X,\delta_*\ca F\big)$
\end{center}

\newpage

\noindent\textit{Proof.} Denote by $\ca O$ the structure presheaf of $X$. By \hyperlink{5.1.2}{Lemma 5.1.2}, \hyperlink{3.2.1.1}{Example 3.2.1.1}, exactness of functor $+$ \cite[3.1]{Toh}, \hyperlink{6.1.1}{Lemma 6.1.1}, \hyperlink{6.1.3}{Proposition 6.1.3} and \hyperlink{6.1.2}{Lemma 6.1.2}, there are natural weak equivalences of complexes of $\ca O_X$-modules
\begin{center}
$\ca C_X=C(\ca O)^+\we\Big(\ca O\overset{\mathbb{L}}{\underset{\ca O^e}{\otimes}} \ca O\Big)^+\we\ca O_X\overset{\mathbb{L}}{\underset{{\ca O_X}^e}{\otimes}} \ca O_X\we\ca O_X\overset{\mathbb{L}}{\underset{\delta^{-1}\ca O_{X\times X}}{\otimes}} \ca O_X\we\mathbb{L}\delta^*\delta_*\ca O_X$
\end{center}
with $C(\ca O):U\map C\big(\ca O(U)\big)$.By exactness, there is a natural weak equivalence $\delta_*\we\mathbb{R\delta_*}$. This allow us to write following natural weak equivalences of chain complexes
\begin{center}
$\mathbb{R}\ba{Hom}_{\ba{Ch}(\ca O_X)}\big(\ca C_X,\ca F\big)\we\mathbb{R}\ba{Hom}_{\ba{Ch}(\ca O_X)}\big(\mathbb{L}\delta^*\delta_*\ca O_X,\ca F\big)$
\end{center} 
\begin{center}
$\we\mathbb{R}\ba{Hom}_{\ba{Ch}(\ca O_X)}\big(\delta_*\ca O_X,\mathbb{R}\delta_*\ca F\big)$
\end{center}
\begin{center}
$\we\mathbb{R}\ba{Hom}_{\ba{Ch}(\ca O_X)}\big(\delta_*\ca O_X,\delta_*\ca F\big)$
\end{center}
\cite[6.7]{Spa}.\\

In fact, \hyperlink{6.1.4}{Theorem 6.1.4} is true as soon as $\delta_*\we\mathbb R\delta_*$. One can see in the next example that it is the obstruction of the veracity of \hyperlink{6.1.4}{Theorem 6.1.4}.\\

\noindent\hypertarget{6.1.5}{\textbf{Example 6.1.5.}} \textit{A non-separated case.} Recall the scheme $X=S\cup_US$ with $S=\ba{Spec}(A)$ and $U=D(s)$ as in \hyperlink{4.0.2}{Example 4.0.2}. We saw that for any $\ca O_X$-module $\ca F$, the cohomology groups $HH^n(X,\ca F$) are isomorphic to the cohomology groups of the derived hom complex $\mathbb{R}\ba{Hom}_{\ba{Ch}(\ca O_S)}\big(\ca C_X,\mathbb Rp_*\ca F\big)$. Consider the following commutative square of schemes
\begin{center}
$\xymatrix @!0 @R=2cm @C=3cm {X \ar[r]^p \ar[d]_{\delta_X} & S \ar[d]^{\delta_S} \\
X\times X \ar[r]_{\Delta p} & S\times S}$
\end{center}
\cite[I, 3.2.1, 5.3.1]{EGA}. There is a natural weak equivalence $\mathbb R(\delta_S)_*\mathbb Rp_*\we\mathbb R(\Delta p)_*\mathbb R(\delta_X)_*$ \cite[6.7]{Spa} and since $S$ is a separated scheme, we have as in \hyperlink{6.1.4}{Theorem 6.1.4} a weak equivalence of complexes of $\ca O_S$-modules $\ca C_S\we\mathbb{L}(\delta_S)^*(\delta_S)_*\ca O_S$. Thus, we have weak equivalences
\begin{center}
$\mathbb{R}\ba{Hom}_{\ba{Ch}(\ca O_S)}\big(\ca C_X,\mathbb Rp_*\ca F\big)\we\mathbb{R}\ba{Hom}_{\ba{Ch}(\ca O_S)}\big(\mathbb{L}(\delta_S)^*(\delta_S)_*\ca O_S,\mathbb Rp_*\ca F\big)$
\end{center}
\begin{center}
$\we\mathbb{R}\ba{Hom}_{\ba{Ch}(\ca O_{S\times S})}\big((\delta_S)_*\ca O_S,\mathbb R(\delta_S)_*\mathbb Rp_*\ca F\big)$
\end{center}
\begin{center}
$\we\mathbb{R}\ba{Hom}_{\ba{Ch}(\ca O_{S\times S})}\big((\delta_S)_*\ca O_S,\mathbb R(\Delta p)_*\mathbb R(\delta_X)_*\ca F\big)$
\end{center}
\cite[6.7]{Spa}. On the other hand, $\ba{Ext}^n_{\ca O_{X\times X}}\big((\delta_X)_*\ca O_X,(\delta_X)_*\ca F\big)$ is the $n$-cohomology of the derived hom complex $\mathbb{R}\ba{Hom}_{\ba{Ch}(\ca O_X)}\big((\delta_X)_*\ca O_X,(\delta_X)_*\ca F\big)$. Since $\delta_S$ is an affine morphism \cite[II, 1.2.1]{EGA} and functor $(\Delta p)^*$ is exact, we have weak equivalences of complexes of $\ca O_{X\times X}$-modules
\begin{center}
$(\delta_X)_*\ca O_X\iso (\delta_X)_*p^*\ca O_S\iso (\Delta p)^*(\delta_S)_*\ca O_S\we\mathbb L(\Delta p)^*(\delta_S)_*\ca O_S$
\end{center}
\cite[7.4]{Swa} and so we have weak equivalences
\begin{center}
$\mathbb{R}\ba{Hom}_{\ba{Ch}(\ca O_X)}\big((\delta_X)_*\ca O_X,(\delta_X)_*\ca F\big)\we\mathbb{R}\ba{Hom}_{\ba{Ch}(\ca O_X)}\big(\mathbb L(\Delta p)^*(\delta_S)_*\ca O_S,(\delta_X)_*\ca F\big)$
\end{center}
\begin{center}
$\we\mathbb{R}\ba{Hom}_{\ba{Ch}(\ca O_{S\times S})}\big((\delta_S)_*\ca O_S,\mathbb R(\Delta p)_*(\delta_X)_*\ca F\big)$
\end{center}
\cite[6.7]{Swa}.\\

\subsection{Spectral sequences}
\
\indent Let $X$ be a scheme with diagonal morphism $\delta:X\arr X\times X$. The composition $\ba{Hom}_{\ca O_{X\times X}}(\delta_*\ca O_X,\ \ \ )=\Gamma\circ\ca Hom_{\ca O_{X\times X}}(\delta_*\ca O_X,\ \ \ )$ gives rise for any $\ca O_X$-module $\ca F$ to the usual \textit{Grothendieck spectral sequence} \cite[2.4.1]{Toh}
\begin{center}
$E^{pq}_2=H^p\big(X\times X,\ca Ext^q_{\ca O_{X\times X}}(\delta_*\ca O_X,\delta_*\ca F)\big)\implies\ba{Ext}^{p+q}_{\ca O_{X\times X}}(\delta_*\ca O_X,\delta_*\ca F)$
\end{center}
Assuming for example that $X$ is a separated smooth scheme, one can show that this spectral sequence is isomorphic to the spectral sequence of \hyperlink{4.1.4}{Proposition 4.1.4}.\\

\noindent\hypertarget{6.2.1}{\textbf{Theorem 6.2.1.}} Let $X$ be a separated scheme with diagonal morphism $\delta:X\arr X\times X$ and $\ca F$ an $\ca O_X$-module. Assume that $\ca Ext^q_{\ca O_X}\big(\ca H_{X,q},\ca F\big)=0$ for any integers $p>0$ and $q\geq 0$. The following spectral sequences are naturally isomorphic
\begin{center}
$\ba{Ext}^p_{\ca O_X}\big(\ca H_{X,q},\ca F\big)\implies HH^{p+q}(X,\ca F)$
\end{center}
\begin{center}
$H^p\big(X\times X,\ca Ext^q_{\ca O_{X\times X}}(\delta_*\ca O_X,\delta_*\ca F)\big)\implies\ba{Ext}^{p+q}_{\ca O_{X\times X}}(\delta_*\ca O_X,\delta_*\ca F)$
\end{center}
\ \\
\textit{Proof.} As in the proof of \hyperlink{4.1.5}{Theorem 4.1.5}, we have by hypothesis a natural isomorphism
\begin{center}
$\ba{Ext}^p_{\ca O_X}\big(\ca H_{X,q},\ca F\big)\iso H^p\Big(X,\ca Hom_{\ca O_X}\big(\ca H_{X,q},\ca F\big)\Big)\iso H^p\Big(X\times X,\delta_*\ca Hom_{\ca O_X}\big(\ca H_{X,q},\ca F\big)\Big)$
\end{center}
for any integers $p>0$ and $q\geq 0$ since $\mathbb{R}\Gamma_{X\times X}\circ\delta_*\we\mathbb{R}\Gamma_{X\times X}\circ\mathbb{R}\delta_*\we\mathbb{R}(\Gamma_{X\times X}\circ\delta_*)=\mathbb{R}\Gamma_X$ \cite[1.3.7]{Hov1}. Moreover, there is a natural isomorphism of $\ca O_X$-modules
\begin{center}
$\ca Hom_{\ca O_X}\big(\ca H_{X,q},\ca F\big)\iso H^q \Big(\mathbb{R}\ca Hom_{\ba{Ch}(\ca O_X)}\big(\ca C_X,\ca F\big)\Big)$
\end{center}
for any integer $q\geq 0$. To see this, choose an injective resolution of $\ca O_X$-modules $\ca F\arr\ca I$. The complex of $\ca O_X$-modules $\ca Hom_{\ba{Ch}(\ca O_X)}\big(\ca C_X,\ca I\big)$ computes the derived hom sheaf $\mathbb{R}\ca Hom_{\ba{Ch}(\ca O_X)}\big(\ca C_X,\ca F\big)$ \cite[6.1]{Spa}. Now we have a double complex of $\ca O_X$-modules $\ca E_0^{pq}=\ca Hom_{\ca O_X}\big(\ca C_{X,q},\ca I^p\big)$ with total cohomology $H\Big(\mathbb{R}\ca Hom_{\ba{Ch}(\ca O_X)}\big(\ca C_X,\ca F\big)\Big)$. Filtering it by columns \cite[5.6.1]{Wei}, the first page $\ca E^1$ is the vertical cohomology which is preserved by exact functors $\ca Hom_{\ca O_X}(\ \ \ ,\ca I^p)$ and the second page $\ca E_2$ is the horizontal cohomology which computes the ext sheaves
\begin{center}
$\ca E_1^{pq}=H^q\Big(\ca Hom_{\ca O_X}\big(\ca C_X,\ca I^p\big)\Big)\iso \ca Hom_{\ca O_X}\big(\ca H_{X,q},\ca I^p\big)$
\end{center}
\begin{center}
$\ca E_2^{pq}\iso H^p\Big(\ca Hom_{\ca O_X}\big(\ca H_{X,q},\ca I\big)\Big)=\ca Ext^p_{\ca O_X}\big(\ca H_{X,q},\ca F\big)\iso{\left\{
    \begin{array}{ll}
        \ca Hom_{\ca O_X}\big(\ca H_{X,q},\ca F\big) & \mbox{if } p=0\\
        0 & \mbox{if } p>0
    \end{array}
\right.}$
\end{center}
This gives the announced isomorphism. Using exactness of functor $\delta_*$ and the proof of \hyperlink{6.1.4}{Theorem 6.1.4}, we can end the proof with natural isomorphisms of $\ca O_{X\times X}$-modules
\begin{center}
$\delta_*\ca Hom_{\ca O_X}\big(\ca H_{X,q},\ca F\big)\iso\delta_*H^q \Big(\mathbb{R}\ca Hom_{\ba{Ch}(\ca O_X)}\big(\ca C_X,\ca F\big)\Big)$
\end{center}
\begin{center}
$\iso H^q \Big(\delta_*\mathbb{R}\ca Hom_{\ba{Ch}(\ca O_X)}\big(\ca C_X,\ca F\big)\Big)$
\end{center}
\begin{center}
$\iso H^q \big(\delta_*\mathbb{R}\ca Hom_{\ba{Ch}(\ca O_X)}(\mathbb{L}\delta^*\delta_*\ca O_X,\ca F)\big)$
\end{center}
\begin{center}
$\iso H^q \big(\mathbb{R}\ca Hom_{\ba{Ch}(\ca O_X)}(\delta_*\ca O_X,\mathbb{R}\delta_*\ca F)\big)$
\end{center}
\begin{center}
$\iso H^q \big(\mathbb{R}\ca Hom_{\ba{Ch}(\ca O_X)}(\delta_*\ca O_X,\delta_*\ca F)\big)=\ca Ext^q_{O_X}(\delta_*\ca O_X,\delta_*\ca F)$
\end{center}
for any integer $q\geq 0$ \cite[6.7]{Spa}.\\

\noindent\hypertarget{6.2.2}{\textbf{Remark 6.2.2.}} Let $X$ be a separated smooth scheme and $\ca F$ an $\ca O_X$-module. For any positive integer $q$, the $\ca O_X$-module $\ca H_{X,q}$ is locally free of finite rank \cite[IV, 17.2.3]{EGAIV4} and so we have $\ca Ext^q_{\ca O_X}\big(\ca H_{X,q},\ca F\big)=0$ for any $p>0$ \cite[4.2.3]{Toh}. Thus we can apply \hyperlink{6.2.1}{Theorem 6.2.1}.\\

\bibliographystyle{alpha}
\bibliography{biblio}

\end{document}